\documentclass[12pt]{article}
\RequirePackage{amsthm, amsmath, amsfonts, amssymb, graphicx}
\RequirePackage[numbers,square]{natbib}
\usepackage{dsfont}
\usepackage{xcolor}

\usepackage{float}
\usepackage[shortlabels]{enumitem}
\usepackage{booktabs}
\newcommand{\bneqn}{\vspace{-0.25cm}\begin{eqnarray}}
\newcommand{\eneqn}{\end{eqnarray}}
\allowdisplaybreaks

\usepackage[utf8]{inputenc}
\usepackage[toc,page,title,titletoc]{appendix}
\usepackage{hyperref,bookmark}
\hypersetup{
  colorlinks=true,
  linkcolor=red,
  citecolor=blue,
  filecolor=blue,
  urlcolor=blue,
}

\newtheorem{theorem}{Theorem}
\newtheorem{lemma}[theorem]{Lemma}

\newtheorem{proposition}{Proposition}
\newtheorem{corollary}{Corollary}

\advance \topmargin by -\headheight
\advance \topmargin by -\headsep
\advance \topmargin .2in
\title{Yule's ``nonsense correlation" for Gaussian random walks}
\author{Philip A. Ernst\footnote{Department of Statistics, Rice University}, Dongzhou Huang\footnote{Department of Statistics, Rice University}, and Frederi G. Viens\footnote{Department of Statistics and Probability, Michigan State University}}
\begin{document}
\maketitle

\begin{abstract}

The purpose of this paper is to provide an exact formula for the second moment of the empirical correlation of two independent Gaussian random walks as well as implicit formulas for higher moments. The proofs are based on a symbolically tractable integro-differential representation formula for the moments of any order in a class of empirical correlations, first established by \cite{ernst2019distribution} and investigated previously  in  \cite{ernst2017yule}. We also provide rates of convergence of the empirical correlation of two independent Gaussian random walks to the empirical correlation of two independent Wiener processes, by exploiting the explicit nature of the computations used for the moments. At the level of distributions, in Wasserstein distance, the convergence rate is the inverse $n^{-1}$ of the number of data points $n$. This holds because we represent and couple the discrete and continuous correlations on a common probability space, where we establish convergence in $L^1$ at the rate $n^{-1}$.
\end{abstract}

\section{Introduction.} \label{sec:introduction}

The main purpose of the present work is to provide an exact formula for the second moment of the empirical correlation of two independent Gaussian random walks, and to apply the method of proof to the question of how fast that correlation converges to that of two independent Wiener processes. We begin this introduction by briefly providing our study's mathematical context.  We then divulge an outline of the mathematical derivation used in our paper for computing the second moment, which is novel and of independent interest. We also provide an outline of the strategy for our convergence result, which is motivated both by numerics presented herein, made possible by our explicit second moment formula, and by the tools developed to compute the second moment.

Let $\{X_{k}\}_{k=1}^{\infty}$ and $\{Y_{k}\}_{k=1}^{\infty}$ be two independent sequences of independent identically distributed random variables with mean $0$ and variance $1$. Define the corresponding partial sums by
\begin{equation}
 S_n = \sum_{j=1}^{n} X_j \quad \text{and} \quad T_n = \sum_{j=1}^{n} Y_j. \label{eq:defforST}
\end{equation}
The empirical correlation of these two random walks is then defined in the usual way as
\begin{equation}
\theta_{n} := \frac{ \frac{1}{n} \sum_{i=1}^{n} S_i T_i - \frac{1}{n^2} (\sum_{i=1}^{n} S_i) (\sum_{i=1}^{n} T_i) }{  \sqrt{ \frac{1}{n} \sum_{i=1}^{n} S_i^2 - \frac{1}{n^2} (\sum_{i=1}^{n} S_i)^2 } \sqrt{ \frac{1}{n} \sum_{i=1}^{n} T_i^2 - \frac{1}{n^2} (\sum_{i=1}^{n} T_i)^2 }}.  \label{eq:defforthetan}
\end{equation}
Despite Udny Yule's warning in \cite{yule1926we} that in the case of two independent
random walks, the observed correlation coefficient has a very different distribution
than that of the nominal $t$-distribution, it has been erroneously assumed that for large enough $n$, these empirical correlations should be small (see \cite{ernst2017yule} and references therein). \\
\indent In \cite{phillips1986understanding}, Phillips  calculated an expression for the limit of these correlations (in the sense of weak convergence), which can be viewed as the empirical correlation of two independent Wiener processes. Henceforth, we shall denote $\theta$ to be the limit of the  correlations $\theta_{n}$. In 2017, Ernst et al. \cite{ernst2017yule} investigated the distribution of the limit $\theta$ by explicitly calculating the standard deviation of the limit to be nearly $0.5$, providing the first formal proof that these correlations $\theta_{n}$ are not small even for arbitrarily large $n$. In 2019, Ernst et al. \cite{ernst2019distribution} succeeded in calculating the moments of $\theta$ up to order 16 and provided the first approximation to the the density of Yule’s ``nonsense correlation''.

In this paper, we explore a question that has not been addressed in the aforementioned references: what is the exact distribution of $\theta_{n}$? This question is of interest not least because discrete stochastic process data (for example, time series data) occur most frequently and extensively in the real world. A test statistic for discrete processes is thus easier for practitioners to apply than that for continuous stochastic processes. Studying the discrete-data test statistic directly is also a means of minimizing the risk of using the continuous statistic abusively when the discrete-data situation is not sufficiently well approximated by a continuous-data one.

However, the task of finding the exact distribution of $\theta_{n}$ for any $n$ and for \textit{any random walk} has proved elusive. Relevant work in this vein includes a series of papers by Andersen (\cite{andersen1953sums,andersen1954fluctuations,andersen1955fluctuations}) provided a combinatorial method based on the idea of cyclic permutations to investigate problems of discrete sequences of partial sums. However, Andersen's methods cannot be applied to evaluate the moments of $\theta_{n}$ since an event generated by $\theta_{n}$ is not invariant under cyclic permutations. The methods used in \cite{phillips1998new} to develop asymptotic theory for spurious regressions, namely, decomposing continuous stochastic process in terms of their orthonormal representations, cannot be employed to find the exact distribution of $\theta_n$ due to the lack of a continuous pattern in the partial sums $S_{k}$ and $T_{k}$.

We now mention a few works which have most directly inspired the present work. In \cite{erd1946certain}, Erd\"{o}s and Kac investigated the asymptotic distributions of four statistics of partial sums of independent identically distributed random variables each having mean $0$ and variance $1$. In \cite{magnus1986exact}, Magnus evaluated the moments of the ratio of a pair of quadratic forms in normal variables, i.e., $x'Ax / x'B x$, where $A$ is symmetric, $B$ positive semidefinite and $x$ is a Gaussian random vector. It is this work in particular which motivates the present work's focus on Gaussian random walks. 

Henceforth, in addition to assuming that $\{X_{k}\}_{k=1}^{\infty}$ and $\{Y_{k}\}_{k=1}^{\infty}$ are two independent sequences of independent random variables, we shall also assume that these variables are standard Gaussian. This specific context will allow us to derive an explicit formula to calculate the second moment of the empirical correlation $\theta_{n}$ for any $n$. 

Our proof of this formula is based on a symbolically tractable integro-differential representation formula for the moments of any order in a class of empirical correlations, established by \cite[Proposition 1]{ernst2019distribution} and investigated previously  in  \cite{ernst2017yule} (see Proposition \ref{propErnst2019} below). The key step in applying this formula is the explicit computation of the joint moment generating function (mgf) $\phi_n$ of the three empirical sums of products and squares which appear in the empirical correlation $\theta_{n}$. This is the topic of Section \ref{sec:jointmgf}. One may also use this representation formula to compute moments of $\theta_{n}$ of any order numerically, using symbolic algebra software. Indeed, we provide these moments up to order 16 for any $n$ via \textsf{Mathematica}. Thus the method for evaluating all moments relies on the joint mgf for the three bilinear and quadratic forms appearing in $\theta_{n}$ (\cite{ernst2019distribution,ernst2017yule}).

The main mathematical contribution in the present paper is the explicit computation of the joint trivariate mgf $\phi_n$ in Section \ref{sec:jointmgf}. To express the second moment of $\theta_{n}$ via the aforementioned representation, it is necessary to compute the partial derivative of $\phi_n$ with respect to its middle variable (the variable representing the empirical covariance). This latter calculation, in Section \ref{sec:momentforguassian}, is only straightforward because of the explicitness of our formula for $\phi_n$. The technical path followed in Section \ref{sec:jointmgf} to compute $\phi_n$ is to express in matrix form the bilinear form mapping the two i.i.d. data sequences $X$ and $Y$ up to the $n$th terms into the empirical covariance of their partial sums $S_n$ and $T_n$, and to compute the matrix's alternative characteristic polynomial $d_n$.  We derive a fully explicit expression for $d_n$ in the Appendix, recursively for $n \ge 5$, by using standard operations to convert $d_n$ into a linear recursion involving a new determinant in tri-diagonal form except for one line along which to expand said new determinant. In doing so, we notice a slight break in the new determinant's recursive nature. When substituting a cell in the determinant's matrix which fixes this break, a second-order recursion emerges, which can be solved explicitly. Relating this back to the original $d_n$ reveals a simple explicit relation, and thus, an explicit formula for $d_n$. 

From a probabilistic standpoint, the Gaussian property of $(X,Y)$ is what allows us to complete this calculation so explicitly. Specifically, we use the properties that the multivariate standard normal law is invariant under orthogonal transformations, and that the Laplace transform of a quadratic form of a bivariate normal vector is a function of a quadratic function. From an analytical standpoint, to compute $d_n$ explicitly, we drew inspiration from the limiting case of $S$ and $T$ distributed as Brownian motions, where Hilbert's approach to Fredholm theory gave us a strong motivation to believe that $d_n$ could be computed. Indeed, the limit of $d_n(\lambda )$ under the appropriate Brownian scaling is explicit, equal to $\sinh (i \sqrt{ \lambda})/(i \sqrt{ \lambda})$ which was a main ingredient in \cite{ernst2017yule}, and also equal, via Mercer's theorem, to $\prod_k (1-\lambda / (k \pi)^2)$ where one recognizes the eigenvalues identified in \cite{ernst2017yule}. 

An ultimate contribution of this paper is the study of the rate of convergence of the empirical correlation $\theta_n$ of Gaussian random walks to the empirical correlation $\theta$ of Wiener processes in Wasserstein distance. Inspired by Hilbert's approach to Fredholm theory, we first construct a ratio $A_{n}/\sqrt{B_n C_n}$ identically distributed with $\theta_n$, where $A_n$, $B_n$ and $C_n$ are second-chaos variables up to constants. We also rewrite $\theta$ as $A/\sqrt{BC}$, where $A$, $B$ and $C$ are also second-chaos variables up to constants. A key element in the setup is to note that, not only can the empirical correlations be represented as ratios involving second-chaos variables, but they can be coupled on the same Wiener space $\Omega$ by using their kernel representations as double integrals with respect to the same pair of independent Wiener processes. Relying on techniques of Wiener chaos, the convergences in $L^2(\Omega)$ of $A_n$, $B_n$ and $C_n$ to $A$, $B$ and $C$ respectively at rate $n^{-2}$ are derived. We then note the Wasserstein distance between $\theta_n$ and $\theta$ is bounded by the $L^1(\Omega)$-norm of $A_{n}/\sqrt{B_n C_n}-A/\sqrt{BC}$, which is
 bounded by a function of the second moments of $A_n -A$, $B_n -B$ and $C_n -C$ and the negative moments of $B_n$, $C_n$, $B$ and $C$. What is left is to give upper bounds for the negative moments. Our idea is to represent these negative moments as a single integral of the product of a positive power function and their moment generating functions (mgfs) and then give upper bounds for mgfs, hence, for negative moments. This idea only works when the mgfs are integrable at $0$ and decay rapidly when approximating to $\infty$. Fortunately, these mgfs follow immediately from the joint mgfs $\phi_n$ and $\phi$ and satisfy the above properties. 
 
 It is worth mentioning that the mgf of $B_n$ or $C_n$ is $1$ over the square root of $d_{n}(-2s/n)$, which is a polynomial with strictly positive coefficients. Furthermore, the coefficients of $d_{n}(-2s/n)$ are eigenvalues of the positive definite matrix $K_n$ (as defined in Section \ref{sec:notations}) after appropriate scaling. We anticipate that these eigenvalues converge to those of the positive definite operator $T_{M}$ defined in \cite{ernst2017yule}. 
 This insight motivates us to establish the existence of a lower bound for $d_{n}(-2s/n)$ for $s\geq 0$ which is uniform for large enough $n$, hence, a uniform upper bound for $E[B_{n}^{-1}]$ and $E[C_{n}^{-1}]$. All of these details are presented in Section \ref{sec:convergenceinWasserstein}. We also conjecture that the entire eigenstructure of $K_n$ converges in some sense to that of $T_M$, though this additional insight is not needed to motivate the uniform lower bound on the negative moments.

The remainder of the paper is organized as follows. In Section \ref{sec:notations}, we introduce necessary notation. In Section \ref{sec:jointmgf}, Theorem \ref{thm1} provides the joint moment generating function needed for obtaining the distribution of $\theta_{n}$ for all $n$. In Section \ref{sec:momentforguassian}, Theorem \ref{thm2} provides an explicit formula for the second moment of $\theta_n$ for any $n$. Numerics for all moments of $\theta_n$ for all $n$ are also given in Section ref{sec:momentforguassian}. The latter motivates our investigation in Section \ref{sec:convergenceinWasserstein} of the rate of convergence of $\theta_n$ to $\theta$. We conclude with Section \ref{sec6}, which provides opportunities for future work which should be tractable given some known tools and techniques in the analysis on Wiener chaos, and could have potential applications to statistical testing based on paths of time series.

\section{Notation.} \label{sec:notations}
We use $I_{n}$ to denote the $n\times n$ identity matrix. For $n \geq 2$ an integer, we define the $(n-1)\times (n-1)$ symmetric matrix $K_{n}$ by 
\begin{equation*}
K_{n} = \left\{ \min(j,k)/n - jk/n^2 \right\}_{j,k=1}^{n-1}, 
\end{equation*}
and its ``alternative characteristic polynomial" $d_{n}(\lambda)$ by 
\begin{equation*}
d_{n}(\lambda) = \det(I_{n-1} -  \lambda K_n). 
\end{equation*} 

We explained in the introduction that the matrix $K_n$ is identified as the discrete-time version of the operator $T_M$ which was identified as a key to the calculations in \cite{ernst2017yule}. In that paper, it was shown that the numerator of the continuous-time Yule's ``nonsense correlation'' $\theta$ (see the definition of $\theta$ in \eqref{theta} in Section \ref{sec:convergenceinWasserstein}) can be written as a member of the second Wiener chaos in its double-Wiener-integral representation, where the bivariate kernel $M$ in that integral is none other than $M(s,t)=\min (s,t)-st$. The expression above for $K_n$ thus comes as no surprise, as the discrete version of $M$. However, as we will see in the next section, $K_n$ also arises naturally when one attempts to express the numerator of $\theta_n$ using the increments $X,Y$ of the random walks $S,T$. That natural phenomenon is exactly the discrete-time analogue of what occurs when identifying the numerator of $\theta$ as a double Wiener integral.

Denoting the eigenvalues of $K_{n}$ generically by $\lambda_2, \cdots, \lambda_n$ (where the numbering starting at 2 is used as a matter of convenience, whose utility will become apparent in the next section), then the alternative characteristic polynomial can be written as
\begin{equation}
 d_{n}(\lambda) =\prod_{j=2}^{n} (1- \lambda_j \lambda). \label{eq:anotherformofacp}  
\end{equation}
We also define two $ (n-1) \times 1$ column random vectors $\mathbf{X}_{n}$ and $\mathbf{Y}_{n}$ by
\begin{equation*}
\mathbf{X}_{n} := \left( X_2, X_3, \cdots, X_n \right)^{\intercal}  \quad \text{and} \quad
\mathbf{Y}_{n} := \left( Y_2, Y_3, \cdots, Y_n \right)^{\intercal}, 
\end{equation*}
where $\{X_{k}\}_{k=1}^{\infty}$ and $\{Y_{k}\}_{k=1}^{\infty}$ are the two independent sequences of independent standard Gaussian random variables used to define the Gaussian random walks $S$ and $T$.
Let
\begin{eqnarray}
&& Z_{11}^n := \frac{1}{n} \sum_{i=1}^{n} S_i^2 - \frac{1}{n^2} \left(\sum_{i=1}^{n} S_i\right)^2, \label{eq:defz11} \\
&& Z_{22}^n := \frac{1}{n} \sum_{i=1}^{n} T_i^2 - \frac{1}{n^2} \left(\sum_{i=1}^{n} T_i\right)^2,  \label{eq:defz22} \\
&& Z_{12}^n := \frac{1}{n} \sum_{i=1}^{n} S_i T_i - \frac{1}{n^2} \left(\sum_{i=1}^{n} S_i\right) \left(\sum_{i=1}^{n} T_i\right), \label{eq:defz12}
\end{eqnarray}
where $S_{i}$ and $T_{i}$ are defined in \eqref{eq:defforST}. Together with \eqref{eq:defforthetan}, we can check easily that 
$$\theta_n = \frac{Z_{12}^n}{\sqrt{Z_{11}^n Z_{22}^n}}.$$ Finally, let us define the joint moment generating function (joint mgf) of the random vector $\left( Z^n_{11}, Z^n_{12}, Z^n_{22} \right)$ by
\begin{equation*}
\phi_{n}(s_{11}, s_{12}, s_{22}) := E \left[ \exp \left\{ -\frac{1}{2} \left( s_{11} Z_{11}^{n} + 2 s_{12} Z_{12}^{n} + s_{22} Z_{22}^{n} \right) \right\} \right], 
\end{equation*}
where $s_{11}, s_{12} $ and $s_{22}$ are such that $s_{11}, s_{22} \geq 0$ and $s_{12}^2 \leq s_{11} s_{22}$.
These inequalities ensure that $\phi_{n}(s_{11}, s_{12}, s_{22})$ is well-defined, as we shall see in Section \ref{sec:jointmgf}. The reader may also check, as a heuristic, that if the possibly ex-centered second-chaos variables $Z_{i,i}^n$ are thought of as independent squares of standard normals, and $Z_{1,2}^n$ is the product of the normals, then condition $s_{12}^2 \leq s_{11} s_{22}$ is necessary.

\section{Calculating the joint moment generating function.} \label{sec:jointmgf}
In this section, we provide an expression for the joint moment generating function $\phi_{n}(s_{11}, s_{12}, s_{22})$. This is the key piece in enabling us to compute the moments of $\theta_{n}$ for all $n$. In the continuous-time setting of \cite{ernst2017yule}, being able to compute this mgf was also a key element, which relied on the fact that the kernel $M(s,t)=\min (s,t)-st$ of the operator $T_M$ was immediately identified as the covariance of the pinned Brownian motion (a.k.a Brownian bridge) on $[0,1]$, for which the eigenvalues happen to be known. In the discrete case herein, there is no such analogous shortcut.

First, by definition, 
\begin{eqnarray}
&& \sum_{i=1}^{n} S_{i} T_{i} = \sum_{i=1}^{n} \left( \sum_{j=1}^{i} X_j \right) \left( \sum_{k=1}^{i} Y_k \right) \notag \\
&=& \sum_{j,k=1}^{n}\, \sum^{n}_{i=\max(j,k)} X_j Y_k = \sum_{j,k=1}^{n} (n- \max(j,k) +1) X_j Y_k.  \label{eq:sumST} 
\end{eqnarray}
Further,
\begin{equation*}
\sum_{i=1}^{n} S_i = \sum_{i=1}^{n} \sum_{j=1}^{i} X_j = \sum_{j=1}^{n} \sum_{i=j}^{n} X_j
= \sum_{j=1}^{n} (n-j+1) X_j.
\end{equation*}
Similarly, 
\begin{equation*}
\sum_{i=1}^{n} T_i = \sum_{k=1}^{n} (n-k+1) Y_k.
\end{equation*}
Hence,
\begin{equation*}
\left(\sum_{i=1}^{n} S_i\right)\left(  \sum_{i=1}^{n} T_i\right)
= \sum_{j,k=1}^{n} (n-j+1)(n-k+1) X_j Y_k.
\end{equation*}
Together with \eqref{eq:defz12} and \eqref{eq:sumST}, we have
\begin{eqnarray}
Z^{n}_{12} &=& \sum_{j,k=1}^{n} \left( \frac{1}{n} \, \Big(n-\max(j,k) +1\Big) - \frac{1}{n^2}\,(n-j+1)(n-k+1) \right) X_j Y_k \notag \\
 &=& \sum_{j,k=1}^{n} \left( \frac{1}{n} \, \Big(\min(j,k) -1\Big) - \frac{1}{n^2} \,(j-1)(k-1) \right) X_j Y_k  \label{eq:quadraticformcoefficients} \\
 &=& \sum_{j,k=2}^{n} \left( \frac{1}{n} \, \Big(\min(j,k) -1\Big) - \frac{1}{n^2} \,(j-1)(k-1) \right) X_j Y_k   \notag \\
 &=& \sum_{j,k=1}^{n-1} \left( \frac{1}{n} \,\min(j,k) - \frac{1}{n^2} \,jk \right) X_{j+1} Y_{k+1}  \notag \\
 &=& \mathbf{X}_{n}^{\intercal} K_{n} \mathbf{Y}_{n},  \notag
\end{eqnarray}
where the third equality holds because $\left(\min(j,k) -1\right)/n - (j-1)(k-1)/n^{2}$ equals to $0$ if either one of the indices $j,k$ is $1$ and the fourth equality holds by making the change of variables $j:=j-1$ and $k := k-1$. As announced in the previous section, we recognize $K_n(j,k)$ defined there and identified here in the last displayed line above, as the discrete version of $M(s,t)=\min(s,t)-st$. Similarly to the expression for $Z_{12}^n$, we have
\begin{equation*}
Z^{n}_{11} = \mathbf{X}_{n}^{\intercal} K_{n} \mathbf{X}_{n}
\quad \text{and} \quad
Z^{n}_{22} = \mathbf{Y}_{n}^{\intercal} K_{n} \mathbf{Y}_{n}. 
\end{equation*}

We now note that since $K_{n}$ is a $(n-1) \times (n-1)$ symmetric matrix, there exits a $(n-1) \times (n-1)$ orthogonal matrix $P_n$ such that
\begin{equation*}
K_{n} = P_{n}^{\intercal} \mathrm{diag} (\lambda_2, \lambda_3, \cdots, \lambda_{n}) P_{n},  
\end{equation*}
where $\lambda_2, \lambda_3, \cdots, \lambda_{n}$ are eigenvalues of $K_n$ and $\mathrm{diag} (\lambda_2, \lambda_3, \cdots, \lambda_{n})$ is a diagonal matrix whose entry in the $j$-th row and the $j$-th column is $\lambda_{j+1}$. Let
\begin{eqnarray*}
 &&\mathbf{\widetilde{X}}_{n} = \left( \widetilde{X}_2, \widetilde{X}_3, \cdots, \widetilde{X}_n \right)^{\intercal}
 := P_{n} \mathbf{X}_{n},  \\
 &&\mathbf{\widetilde{Y}}_{n} = \left( \widetilde{Y}_2, \widetilde{Y}_3, \cdots, \widetilde{Y}_n \right)^{\intercal}
 := P_{n} \mathbf{Y}_{n},  
\end{eqnarray*}
be two $(n-1)\times 1$ column random vectors. Since $\mathbf{X}_{n}$ and $\mathbf{Y}_{n}$ are two independent Gaussian random vectors with distribution $\mathcal{N}\left( \mathbf{0}, I_{n-1} \right)$ and because $P_n$ is an orthogonal matrix, then $\mathbf{\widetilde{X}}_{n}$ and $\mathbf{\widetilde{Y}}_{n}$ are also two independent Gaussian random vectors with distribution $\mathcal{N}\left( \mathbf{0}, I_{n-1} \right)$. This implies that $\widetilde{X}_2, \widetilde{X}_3, \cdots, \widetilde{X}_n$, $\widetilde{Y}_2, \widetilde{Y}_3, \cdots, \widetilde{Y}_n$ are independent standard Gaussian random variables.

Before presenting our formula for the trivariate mgf $\phi_n$ in Theorem \ref{thm1} below, we reveal an explicit calculation of the alternative characteristic polynomial $d_{n}(\lambda)$. The proof is relegated to the Appendix.

\begin{lemma}\label{lem1}
The alternative characteristic polynomial $d_{n}(\lambda)$ may be written as
\begin{eqnarray}
d_{n}(\lambda)
&=&  \frac{1}{n \sqrt{\left(\frac{\lambda}{n}-2\right)^2 -4}}\left( - \dfrac{\left(\frac{\lambda}{n} -2\right) - \sqrt{\left(\frac{\lambda}{n}-2\right)^2 - 4}}{2} \right)^{n}  \notag  \\
&& -\frac{1}{n \sqrt{\left(\frac{\lambda}{n}-2\right)^2 -4}}
  \left( - \frac{\left(\frac{\lambda}{n} -2\right) + \sqrt{\left(\frac{\lambda}{n}-2\right)^2 - 4}}{2} \right)^{n}  \notag \\ [1mm]
&=& \frac{(-1)^{n-1}}{n \cdot 2^{n-1} } \, \sum_{k=1}^{\lceil n/2 \rceil} {n \choose 2k-1} \, \left(\frac{\lambda}{n}-2\right)^{n-(2k-1)} \, \left( \left(\frac{\lambda}{n}-2\right)^2 -4 \right)^{k-1}, \label{eq:representationfordinmaindocument}
\end{eqnarray}
where $\lceil x \rceil$ is the least integer greater than or equal to $x$.
\end{lemma}
\begin{proof}
See Appendix.
\end{proof}

With above in hand, we now calculate the joint mgf $\phi_{n}$.

\begin{theorem}\label{thm1}
The joint moment generating function $\phi_n$ for the triple $(Z_{11}^n,Z_{12}^n,Z_{22}^n)$ of random variables defined in \eqref{eq:defz11}, \eqref{eq:defz12}, \eqref{eq:defz22} is given for $s_{11}, s_{22} \geq 0$ and $s_{12}^2 \leq s_{11} s_{22}$, by
\begin{equation*}
\phi_{n}(s_{11}, s_{12}, s_{22})= \left( d_{n}(\alpha) \, d_{n}(\beta)  \right)^{-1/2} 
\end{equation*}
where $\alpha$ and $\beta$ are defined as follows:
\begin{eqnarray}
&&\alpha := \alpha\left(s_{11}, s_{12}, s_{22}\right) = -\frac{s_{11}+s_{22} + \sqrt{(s_{11}-s_{22})^2 + 4 s_{12}^{2}}}{2} ,  \label{eq:alpha} \\ [1mm]
&&\beta :=  \beta\left(s_{11}, s_{12}, s_{22}\right) = -\frac{s_{11}+s_{22} - \sqrt{(s_{11}-s_{22})^2 + 4 s_{12}^{2}}}{2},  \label{eq:beta} 
\end{eqnarray}
\end{theorem}

\begin{proof}
We first calculate 
\begin{eqnarray*}
&&s_{11} Z^{n}_{11} + 2 s_{12} Z^{n}_{12} +  s_{22} Z^{n}_{22}  \\ [1mm]
&=& s_{11} \mathbf{X}_{n}^{\intercal} K_{n} \mathbf{X}_{n} + 2 s_{12} \mathbf{X}_{n}^{\intercal} K_{n} \mathbf{Y}_{n} +  s_{22} \mathbf{Y}_{n}^{\intercal} K_{n} \mathbf{Y}_{n}   \\ [1mm]
&=& s_{11}\, \mathbf{\widetilde{X}}_{n}^{\intercal}\, \mathrm{diag} (\lambda_2, \lambda_3, \cdots, \lambda_{n}) \, \mathbf{\widetilde{X}}_{n}
+ 2s_{12}\, \mathbf{\widetilde{X}}_{n}^{\intercal}\, \mathrm{diag} (\lambda_2, \lambda_3, \cdots, \lambda_{n}) \, \mathbf{\widetilde{Y}}_{n} \\  [1mm]
&& + s_{22}\, \mathbf{\widetilde{Y}}_{n}^{\intercal}\, \mathrm{diag} (\lambda_2, \lambda_3, \cdots, \lambda_{n}) \, \mathbf{\widetilde{Y}}_{n}  \\
&=& s_{11} \,\sum_{j=2}^{n} \lambda_{j} \widetilde{X}_{j}^{2}  + 2s_{12} \,\sum_{j=2}^{n} \lambda_{j} \widetilde{X}_{j} \widetilde{Y}_{j} + s_{22} \,\sum_{j=2}^{n} \lambda_{j} \widetilde{Y}_{j}^{2}  \\
&=& \sum_{j=2}^{n} \lambda_{j} \left( s_{11} \widetilde{X}_{j}^{2} + 2 s_{12} \widetilde{X}_{j} \widetilde{Y}_{j} + s_{22} \widetilde{Y}_{j}^{2} \right). 
\end{eqnarray*}
By independence of $\widetilde{X}_{j}$ and $\widetilde{Y}_{k}$ for $j,k \in \{2,3,\cdots,n\}$,
\begin{eqnarray}
\phi_{n}(s_{11}, s_{12}, s_{22}) &=&  E \left[ \exp \left\{ -\frac{1}{2} \left( s_{11} Z_{11}^{n} + 2 s_{12} Z_{12}^{n} + s_{22} Z_{22}^{n} \right) \right\} \right]  \notag \\
&=& \prod_{j=2}^{n} E \left[ \exp\left\{ -\frac{1}{2} \lambda_{j} \left( s_{11} \widetilde{X}_{j}^{2} + 2 s_{12} \widetilde{X}_{j} \widetilde{Y}_{j} +  s_{22} \widetilde{Y}_{j}^{2} \right) \right\} \right] \notag \\
&=& \prod_{j=2}^{n} \left( 1 + (s_{11} + s_{22}) \lambda_{j} + (s_{11}s_{22} - s_{12}^{2}) \lambda_{j}^{2} \right)^{-1/2}  \label{eq:laplaceofquadraticform}  \\
&=& \prod_{j=2}^{n} \left( (1- \alpha \lambda_j) (1- \beta \lambda_j) \right)^{-1/2} \notag \\
&=& \left( \prod_{j=2}^{n} (1-\alpha \lambda_j) \,  \prod_{j=2}^{n} (1-\beta \lambda_j) \right)^{-1/2}  \notag \\
&=& \left( d_{n}(\alpha) \, d_{n}(\beta)  \right)^{-1/2}, \label{eq:phibyd}  
\end{eqnarray}
where $\alpha$ and $\beta$ are defined in the statement of the theorem. Note that in line \eqref{eq:laplaceofquadraticform} a standard expression for the mgf of a linear-quadratic functional of a normal variable has been used (iteratively twice), and, further, the independence of $\widetilde{X}_j$ and $\widetilde{Y}_j$ has been employed. Note that in line \eqref{eq:laplaceofquadraticform} the conditions $s_{11}, s_{22} \geq 0$ and $s_{12}^{2} \leq s_{11} s_{22}$ ensure the applicability of the standard expression for the mgf of a linear-quadratic functional of a bivariate random vector. The last equality holds by the representation of the alternative characteristic polynomial of $K_{n}$ by the eigenvalues of $K_n$, see \eqref{eq:anotherformofacp}. Combining \eqref{eq:representationfordinmaindocument}, \eqref{eq:alpha}, and \eqref{eq:beta} allows us to represent the joint mgf $\phi_{n}(s_{11}, s_{12}, s_{22})$ explicitly in terms of $d_{n}(\lambda)$, $\alpha(s_{11}, s_{12}, s_{22})$ and $\beta(s_{11}, s_{12}, s_{22})$, as is given by \eqref{eq:phibyd}, and as announced in Theorem \ref{thm1}.
\end{proof}

\section{Moments of $\theta_{n}$.} \label{sec:momentforguassian}
In the previous section, we gave an exact representation for the joint trivariate mgf $\phi_{n}$. In this section, we use it to calculate the moments of $\theta_{n}$ by a method provided by Ernst et al. (see Proposition 1 in \cite{ernst2019distribution}), which we cite as follows:
\begin{proposition}[Ernst et al. (2019)]  \label{propErnst2019}
For $m = 0,1,2,\cdots$, we have
\begin{equation}
E \left( \theta_{n}^{m} \right) = \frac{(-1)^{m}}{2^{m} \Gamma(m/2)^{2} }
\int_{0}^{\infty} \int_{0}^{\infty} s_{11}^{m/2-1} s_{22}^{m/2-1}\, \frac{\partial^{m}\phi_{n}}{\partial s_{12}^{m}}(s_{11}, 0 , s_{22})\, ds_{11} ds_{22}. \label{eq:allmoments}  
\end{equation}
\end{proposition}
\bigskip
An immediate application of this proposition yields that the second moment of $\theta_{n}$ is given by the following double Riemann integral:
\begin{equation}
E \left( \theta_{n}^{2} \right) = \frac{1}{4} \int_{0}^{\infty} \int_{0}^{\infty} \frac{\partial^2 \phi_{n}}{\partial s_{12}^{2}}(s_{11}, 0,  s_{22}) \, ds_{11} ds_{22}. \label{eq:thetasecondmoment} 
\end{equation}

\subsection{Explicit formula for the second moments of $\theta_{n}$.}
We now calculate the integrand in the previous integral representation explicitly, yielding the next theorem, which is a closed-form expression for the second moment of $\theta_n$ for any $n$. 

\begin{theorem}\label{thm2}
The second moment of $\theta_n$ is
\begin{eqnarray}
&&E \left( \theta_{n}^{2} \right)  \notag \\
&=& - \frac{1}{4} \int_{0}^{\infty} \int_{0}^{\infty} \frac{(s_{11}+2n)(s_{22} + 2n) + 4 n^2}{\left[ s_{11} s_{22} (s_{11}+4n) (s_{22}+4n) \right]^{3/4}} \left[ f(s_{11}/n)^{n} - f(s_{11}/n)^{-n} \right]^{-1/2} \notag \\
&& \quad \quad \quad \times \left[ f(s_{22}/n)^{n} - f(s_{22}/n)^{-n} \right]^{-1/2}\, ds_{11} ds_{22}  \notag \\ [1mm]
&& + \frac{1}{4} \int_{0}^{\infty} \int_{0}^{\infty} \, \Bigg\{ \frac{n(s_{11} + s_{22} + 4n) }{\sqrt{s_{11}^{2} + 4n s_{11}} + \sqrt{s_{22}^{2} + 4ns_{22}}}  \cdot \big[ f(s_{11}/n)^{n} f(s_{22}/n)^{n} \notag \\
&& \quad \quad - f(s_{11}/n)^{-n} f(s_{22}/n)^{-n} \big]
  + \frac{1}{2} \bigg( \sqrt{s_{11}^{2} + 4n s_{11}} + \sqrt{s_{22}^{2} + 4ns_{22}} + s_{11} + s_{22} +4n \bigg) \notag \\
&& \quad \quad\cdot \frac{  f(s_{11}/n)^{n} f(s_{22}/n)^{-n} - f(s_{22}/n)^{n} f(s_{11}/n)^{-n}  }{f(s_{11}/n)- f(s_{22}/n)} \Bigg\} \cdot \left[ f(s_{11}/n)^{n} - f(s_{11}/n)^{-n} \right]^{-3/2} \notag \\
&& \quad \quad \cdot \left[ f(s_{22}/n)^{n} - f(s_{22}/n)^{-n} \right]^{-3/2} \cdot \left[ s_{11} s_{22} (s_{11} +4n)(s_{22}+4n) \right]^{-1/4} \, ds_{11} ds_{22},  \label{eq:exactsecondmomentoftheta}  
\end{eqnarray}
where 
\begin{equation}\label{eqf}
f(\lambda) := \frac{(\lambda+2) + \sqrt{(\lambda+2)^2 -4} }{2}. 
\end{equation}
\end{theorem}
\begin{proof}
It is sufficient to provide the announced closed form expression for $\frac{\partial^2 \phi_{n}}{\partial s_{12}^{2}}(s_{11}, 0,  s_{22})$. Recalling the definition of $f(\lambda)$ in \eqref{eqf}, a direct calculation yields
\begin{equation}
d_{n}(-\lambda) = \frac{1}{n \sqrt{(\lambda/n +2)^2 -4}} \left[ f(\lambda/n)^{n} - f(\lambda/n)^{-n} \right],  \label{eq:dbyf} 
\end{equation}
\begin{eqnarray}
d_{n}'(-\lambda) &=& \frac{1}{n^2} \, \frac{\lambda/n +2}{\left[(\lambda/n+2)^2 -4\right]^{3/2} } \left[ f(\lambda/n)^{n} - f(\lambda/n)^{-n} \right]  \notag \\
&& - \frac{1}{n} \, \frac{1}{(\lambda/n+2)^2 -4} \left[ f(\lambda/n)^{n} + f(\lambda/n)^{-n} \right],  \label{eq:dprimebyf}
\end{eqnarray}
and
\begin{eqnarray}
\frac{\partial^2 \phi_{n}}{\partial s_{12}^2} &=&
\frac{3}{4} \big( d_{n}(\alpha) d_{n}(\beta)  \big)^{-5/2} \, \left(d_{n}'(\alpha) \, d_{n}(\beta) \, \frac{\partial \alpha}{\partial s_{12}} + d_{n}(\alpha) \, d_{n}'(\beta) \, \frac{\partial \beta}{\partial s_{12}} \right)^{2} \notag \\
&& - \frac{1}{2} \big( d_{n}(\alpha) d_{n}(\beta)  \big)^{-3/2} \,
  \sum_{j=0}^{2} {2 \choose j} \,  d_{n}^{(j)}(\alpha) \, d_{n}^{(2-j)}(\beta) \left(\frac{\partial \alpha}{\partial s_{12}}\right)^{j} \left(\frac{\partial \beta}{\partial s_{12}}\right)^{2-j}  \notag \\
&& - \frac{1}{2} \big( d_{n}(\alpha) d_{n}(\beta)  \big)^{-3/2} \,
  \left(d_{n}'(\alpha) \, d_{n}(\beta) \, \frac{\partial^2 \alpha}{\partial s_{12}^{2}} + d_{n}(\alpha) \, d_{n}'(\beta) \, \frac{\partial^2 \beta}{\partial s_{12}^{2}} \right), \label{eq:phisecondderivative}  
\end{eqnarray}
where $\alpha$ and $\beta$ are defined in \eqref{eq:alpha} and \eqref{eq:beta}. Note that
\begin{eqnarray*}
&& - \frac{\partial \alpha}{\partial s_{12}} = \frac{\partial \beta}{\partial s_{12}} = \frac{2 s_{12}}{\sqrt{(s_{11}-s_{22})^2 + 4 s_{12}^2}}  \\
&& - \frac{\partial^2 \alpha}{\partial s_{12}^{2}} = \frac{\partial^2 \beta}{\partial s_{12}^{2}}
= \frac{2(s_{11}-s_{22})^2}{\left[ (s_{11}-s_{22})^2 + 4s_{12}^{2} \right]^{3/2}}. 
\end{eqnarray*}
It follows easily that $\alpha(s_{11}, 0 , s_{22}) = - \max(s_{11}, s_{22})$, $\beta(s_{11}, 0 , s_{22}) = - \min(s_{11}, s_{22})$, and
\begin{eqnarray*}
&& - \frac{\partial \alpha}{\partial s_{12}} (s_{11}, 0 , s_{22}) = \frac{\partial \beta}{\partial s_{12}} (s_{11}, 0 , s_{22}) =0  \\
&& - \frac{\partial^2 \alpha}{\partial s_{12}^{2}} (s_{11}, 0 , s_{22}) = \frac{\partial^2 \beta}{\partial s_{12}^{2}} (s_{11}, 0 , s_{22}) = 2 |s_{11} - s_{22}|^{-1}. 
\end{eqnarray*}
Plugging the above results into \eqref{eq:phisecondderivative} yields
\begin{eqnarray}
\frac{\partial^2 \phi_{n}}{\partial s_{12}^2}(s_{11}, 0, s_{22})
= \frac{ \frac{d_{n}'(-\max(s_{11}, s_{22}))}{ d_{n} (-\max(s_{11}, s_{22}))} -   \frac{d_{n}'(-\min(s_{11}, s_{22}))}{ d_{n} (-\min(s_{11}, s_{22}))}}{\left[ d_{n}(-\max(s_{11}, s_{22})) d_{n}(-\min(s_{11}, s_{22})) \right]^{1/2} |s_{11}-s_{22}|}.  \label{eq:phi2plugvalue}
\end{eqnarray}
Combining \eqref{eq:thetasecondmoment}, \eqref{eq:dbyf}, \eqref{eq:dprimebyf} and \eqref{eq:phi2plugvalue} and performing straightforward calculations, we arrive at an explicit formula for  $\frac{\partial^2 \phi_{n}}{\partial s_{12}^{2}}(s_{11}, 0,  s_{22})$ and the announced double integral expression for the second moment of $\theta_n$ for all $n$.
\end{proof}

\subsection{Numerics.}
We now turn to numerics. \textsf{Mathematica} allows us to calculate the second moment of Yule's ``nonsense correlation'' $\theta_{n}$ for any given $n$. The numerical results are summarized in Table \ref{tab:Yulenumerical}.
\begin{table}[H]
\centering
\begin{tabular}{ccccccc}
\toprule[1.5pt]
 $n$  & 2  &   5  & 10  &  20  &  50  &  100  \\
 $E \left(\theta_{n}^{2} \right)$  & 1.000000  & 0.341109 & 0.265140 & 0.246645 & 0.241501 & 0.240767 \\
 \hline
 $n$  & 200  &   500  & 1000  &  2000  &  5000  &  $\infty$  \\
 $E \left(\theta_{n}^{2} \right)$  & 0.240584 & 0.240532 & 0.240525 & 0.240523 & 0.240523 & 0.240523 \\
\bottomrule[1.5pt]
\end{tabular}
\caption{Numerical Results of the second moment of $\theta_n$ for various values of $n$}
\label{tab:Yulenumerical}
\end{table}
For higher-order moments as represented in \eqref{eq:allmoments}, we can use \textsf{Mathematica} to perform symbolic high-order differentiation and then the two dimensional integration, thereby implicitly calculating higher moments of $\theta_n$ for all $n$. The numerical results of some higher-order moments of $\theta_{50}$ are summarized in Table \ref{tab:Yulenumericalhighmomentfor50}.
\begin{table}[H]
\centering
\begin{tabular}{ccccc}
\toprule[1.5pt]
 $k$  & 2  &   4  & 6  &  8    \\
 $E \left(\theta_{50}^{k} \right)$  & 0.241501  & 0.109961 & 0.061465 & 0.038257  \\
 \hline
 $k$  & 10  &   12  & 14  &  16   \\
 $E \left(\theta_{50}^{k} \right)$  & 0.025485 & 0.017803 & 0.012885 & 0.009586  \\
\bottomrule[1.5pt]
\end{tabular}
\caption{Numerical Results of higher-order moments of $\theta_{50}$}
\label{tab:Yulenumericalhighmomentfor50}
\end{table}

\section{Convergence in Wasserstein distance.} \label{sec:convergenceinWasserstein}

Table \ref{tab:Yulenumerical} and Table \ref{tab:Yulenumericalhighmomentfor50} in the last section give us insight into the behavior of the distribution of $\theta_n$ for large $n$, as it approximates the distribution of its limit $\theta$ defined below in \eqref{theta}. In Table \ref{tab:Yulenumerical}, we note the rather rapid convergence of  $E \left(\theta_{n}^{2} \right)$ as $n \to \infty$. We observe that this convergence occurs at a rate which appears to be faster than $n^{-1}$. This fact motivates us to investigate the rate of convergence of $\theta_n$ to $\theta$. In this section, we give an upper bound for the Wasserstein distance between $\theta_n$ and $\theta$, which comes from a coupling of $\theta_n$ and $\theta$ on the same probability space $\Omega$, in which the convergence occurs in $L^1(\Omega)$.

Let $W_1$ and $W_2$ be two independent Wiener processes. Then Yule's ``nonsense correlation'' (see \cite{ernst2017yule}) is given by 
\begin{equation}
\label{theta}
\theta = \frac{\int_{0}^{1} W_{1}(t)W_{2}(t) dt - \int_{0}^{1} W_{1}(t) dt \int_{0}^{1} W_{2}(t) dt }{\sqrt{\int_{0}^{1} W_{1}^{2}(t) dt  - \left(\int_{0}^{1} W_{1}(t) dt \right)^{2} } \sqrt{\int_{0}^{1} W_{2}^{2}(t) dt  - \left(\int_{0}^{1} W_{2}(t) dt \right)^{2} } }.
\end{equation}

If $X,Y$ are two real-valued random variables, recall that the Wasserstein distance between the law of $X$ and the law of $Y$ is given by
\begin{equation*}
d_{W} (X,Y) := \sup_{f \in \mathrm{Lip}(1)} |Ef(X) - Ef(Y)|, 
\end{equation*}
where $\mathrm{Lip}(1)$ is the set of all Lipschitz functions with Lipschitz constant $\leq 1$. Our key result (Theorem \ref{thm:boundwassertain}) regarding the convergence of $\theta_n$ to $\theta$ is as follows:
\begin{equation}
d_{W}(\theta_n, \theta) = \mathcal{O}\left(\frac{1}{n}\right). 
\end{equation}
We prove this by showing that $E\left[ \left|\theta_n - \theta \right| \right] = \mathcal{O}\left(\frac{1}{n}\right)$ under a natural coupling of $\theta_n$ and $\theta$. 

The reader will find some heuristic comments regarding how this result arises, and what more could be expected for other processes, at the end of the next, short, subsection which provides the preparatory setup needed to prove Theorem \ref{thm:boundwassertain}.

\subsection{Notation, coupling, extensions and implications}\label{sec:coupling}

Define $M(s,t) := \min(s,t) - st$. For every $n \in \mathds{N}_{+}$, define
\begin{equation}\label{emenn}
M_{n}(s,t) := \sum_{1\leq j,k \leq n} M\left(\frac{j-1}{n}, \frac{k-1}{n}\right) \mathds{1}_{\{(j-1)/n < s \leq j/n\}}
\mathds{1}_{\{(k-1)/n < t \leq k/n\}}. 
\end{equation}
For every $n \in \mathds{N}_{+}$, define
\begin{eqnarray*}
A_n &=& \int_{0}^{1} \int_{0}^{t} M_{n}(s,t) \, dW_{1}(s) \, dW_{2}(t)
+ \int_{0}^{1} \int_{0}^{s} M_{n}(s,t) \, dW_{2}(t) \, dW_{1}(s), \\
B_{n} &=& 2\int_{0}^{1} \int_{0}^{t} M_{n}(s,t) \, dW_{1}(s) \, dW_{1}(t)
 + \frac{1}{n} \sum_{j=1}^{n} M\left(\frac{j-1}{n}, \frac{j-1}{n}\right), \\
C_{n} &=& 2\int_{0}^{1} \int_{0}^{t} M_{n}(s,t) \, dW_{2}(s) \, dW_{2}(t)
+ \frac{1}{n} \sum_{j=1}^{n} M\left(\frac{j-1}{n}, \frac{j-1}{n}\right). 
\end{eqnarray*}
Further, let
\begin{eqnarray*}
A &=& \int_{0}^{1} \int_{0}^{t} M(s,t) \, dW_{1}(s) \, dW_{2}(t)
+ \int_{0}^{1} \int_{0}^{s} M(s,t) \, dW_{2}(t) \, dW_{1}(s), \\
B &=& 2\int_{0}^{1} \int_{0}^{t} M(s,t) \, dW_{1}(s) \, dW_{1}(t)
+  \int_{0}^{1} M(t,t) \, dt , \\
C &=& 2 \int_{0}^{1} \int_{0}^{t} M(s,t) \, dW_{2}(s) \, dW_{2}(t)
+   \int_{0}^{1} M(t,t) \, dt . 
\end{eqnarray*}
A key point here, mentioned in the introduction, is that we choose to use the same pair $(W_1,W_2)$ of independent Wiener processes to represent all six of these variables. This is a natural coupling on the common Wiener space $\Omega$ defined by this pair, which allows us to relate the two empirical correlations to each other in a way that easily yields the Wasserstein distance between their distributions. In particular, in Section \ref{sec:propertyofAnBnCnABC}, we will show $\theta_{n} \overset{D}{=}A_n/\sqrt{B_n C_n}$ and $\theta = A/\sqrt{B C}$, while in the first step in the proof of Theorem \ref{thm:boundwassertain} in Section \ref{sec:UBWass}, we establish that $d_{W}(\mathcal{L}(X),\mathcal{L}(Y)) \leq E[|X-Y|]$ for any pair of integrables rv's $(X,Y)$ on the same probability space, from which we conclude
\begin{equation*}
d_{W}(\theta_n, \theta) = d_{W} \left( \mathcal{L}\left(\frac{A_n}{\sqrt{B_n C_n}} \right),  \mathcal{L}\left(\frac{A}{\sqrt{B C}} \right) \right) \leq E\left[ \left| \frac{A_n}{\sqrt{B_n C_n}} -  \frac{A}{\sqrt{B C}} \right| \right]. 
\end{equation*}
In the sequel, we shall restrict our attention to $A_n, B_n, C_n, A, B$ and $C$. One key reason for defining $A_n, B_n, C_n, A, B$ and $C$ is that, being defined as second-chaos variable plus a constant (which may be 0), the upper bounds for the second moments of $A_n - A$, $B_n -B$ and $C_n -C$ can be estimated as $\mathcal{O}(n^{-2})$. Hence, $A_n$, $B_n$ and $C_n$ converge in $L^2(\Omega)$ to $A$, $B$ and $C$ respectively at rate $\sqrt{n^{-2}}$. In fact, the second moments of $A_n - A$, $B_n -B$ and $C_n -C$ can be calculated explicitly. These details will be stated and proved in Section \ref{sec:propertyofAnBnCnABC}.

This convergence rate $\mathcal{O}(n^{-2})$ converts into the rate $\mathcal{O}(n^{-1})$ in Theorem \ref{thm:boundwassertain} because of the need to separate numerator from denominator. One might view this as the cost to pay for this conversion. However, we think it is more fruitful to view the rate of convergence at the level of norms, which preserve scales: the $\mathcal{O}(n^{-1})$ is the rate of convergence of the three elements constituting $\theta_n$ in $L^2(\Omega)$-norm. This leads us to presume this Wasserstein-distance rate of convergence is sharp, though it is beyond the scope of this paper to establish this rigorously. From the so-called property of hyper-contractivity on fixed Wiener chaos (see \cite{NourdinPeccati} Chapter 2), for all $p > 1$, all $L^p(\Omega)$-norms of the three differences $A_n-A, B_n-B, C_n-C$ are equivalent, making it unnecessary to speculate whether computing the convergence rates of any specific higher moment might provide additional insight. Expanding the ratios defining $\theta_n - \theta$ into tri-variate Taylor series did not lead us to any further insight based on those explicit norm-equivalence universal constants. 

In our mind, it is more interesting to ask oneself whether the rate $\mathcal{O}(n^{-1})$ for $d_W(\theta_n,\theta)$, which is inherited from the rate of $\mathcal{O}(n^{-2})$ for $Var(A_n-A),\,Var(B_n-B),$ and $Var(C_n-C)$, is generic, or whether it is specific to random walks. Resolving this question rigorously is also beyond the scope of this paper, but our preliminary calculations indicate that the aforementioned $\mathcal{O}(n^{-2})$ only holds because of the property of independence of increments of a random walk defined as a partial sum of a sequence of independent terms. 

We believe that for other Markov chains which might converge in law, and specifically for any reasonable discretization of processes which are far from having independent increments, such as long-memory processes or mean-reverting processes, the rate of convergence to 0 of $Var(A_n-A),\,Var(B_n-B),$ and $Var(C_n-C)$ is $\mathcal{O}(n^{-1})$. Using a simple polarization argument, the rates of convergence to 0 of these three differences should be essentially equivalent, so that looking at merely one of them would give the order for all of them. As mentioned elsewhere (e.g. Section \ref{sec6}), the method of proof below to establish the rate $\mathcal{O}(n^{-2})$ is one of direct calculations, but the same rate can also be established using a more generic, less precise calculation where one compares Riemann integrals to their approximations using step functions. When attempting that calculation, the property of independence of increments comes plainly into view, implying a number of cancellations much like what one observes when computing the quadratic variation of a martingale. This same methodology seems to indicate that no such cancellations occur for non-independent-increment cases, but that our conjectured rate $\mathcal{O}(n^{-1})$ is straightforward to establish for other Gaussian processes, using the same type of coupling as for Gaussian random walks. Extending the conjecture to non-Gaussian processes would require more work.

The distinction which we conjecture above between Gaussian random walks and other Gaussian processes is far from merely academic. It means that the use of the properties of the continuous-time limit of Yule's ``nonsense correlation'' $\theta$, which are straightforward to establish using simulations, to infer statistical properties of discrete-time random walks, is legitimate for moderate sample sizes, but not so if the data does not behave like the path of a random walk with independent increments. For instance, a statistic on $\theta$ that relates to the construction of the Wasserstein metric (e.g. a mean value or a moment) can be presumed generically accurate at a 1\% level for a Gaussian random walk with several hundred data points, while tens of thousands of data points would be needed, according to our conjecture, when working with a mean-reverting time series. In the environmental sciences, where such time series are ubiquitous, and where many have yearly frequencies, no such reliance on $\theta$ directly can be assumed on a historical scale. In other application domains, such as in quantitative finance, high-frequency studies over several years, such as when studying the long-term distribution and movements of interest rates or of market volatility, can routinely draw on enough data points, however. Also in finance, shorter-term studies of other objects, like stock returns, relate more readily to Gaussian random walks, where our results herein indicate that only hundreds of measurements over time would allow the use of $\theta$'s law instead of needing to rely on $\theta_n$. For random-walk time series which are shorter yet, our explicit results on $\theta_n$ from Section \ref{sec:momentforguassian} are available.

\subsection{Properties of $A_n, B_n, C_n, A, B$, and $C$.} \label{sec:propertyofAnBnCnABC}
In this section, we derive several properties of $A_n, B_n, C_n, A, B$, and $C$, including justifying the coupling, some convenient a.s. constraints we point out, explicit formulae for univariate moment-generating functions, and most importantly from the standpoint of our quantitative analysis, the last two propositions in this section provide the aforementioned convergences to zero of the variances of the differences between the approximating and limiting three elements constituting $\theta_n$ and $\theta$. 
\begin{proposition} \label{prop:identicaldistribution}
The following statements hold, where the equality in (a) and the first equality in (c) are in distribution:
\begin{enumerate}[(a)]
\item $(Z_{11}^{n}/n, Z_{12}^{n}/n, Z_{22}^{n}/n )  \overset{D}{=} (B_n, A_n, C_n)$ for every $n \in \mathds{N}_{+}$;
\item
\begin{eqnarray}
A &=& \int_{0}^{1} W_{1}(t)W_{2}(t) dt - \int_{0}^{1} W_{1}(t) dt \int_{0}^{1} W_{2}(t) dt , \label{eq:A} \\
B &=& \int_{0}^{1} W_{1}^{2}(t) dt  - \left(\int_{0}^{1} W_{1}(t) dt \right)^{2} , \label{eq:B} \\
C &=& \int_{0}^{1} W_{2}^{2}(t) dt  - \left(\int_{0}^{1} W_{2}(t) dt \right)^{2}. \label{eq:C}
\end{eqnarray}
\item $\theta_{n} \overset{D}{=}A_n/\sqrt{B_n C_n}$ and $\theta = A/\sqrt{B C}$.
\end{enumerate}
\end{proposition}
\begin{proof}
See Appendix.
\end{proof}

A helpful corollary of Proposition \ref{prop:identicaldistribution} is as follows. 
\begin{corollary} \label{coro:1}
(a) $| A_n/ \sqrt{B_n C_n} | \leq 1$ a.s.;
(b) $| A/ \sqrt{B C} | \leq 1$ a.s.;
(c) $B_n, C_n >0$ a.s. for $n \geq 2$;
(d) $B,C >0$ a.s..
\end{corollary}
\begin{proof}
By the Cauchy-Schwarz inequality, $|Z_{12}^{n}/\sqrt{Z_{11}^{n} Z_{22}^{n}}| \leq 1$. By Proposition \ref{prop:identicaldistribution}, $(Z_{11}^{n}/n, Z_{12}^{n}/n, Z_{22}^{n}/n )  \overset{D}{=} (B_n, A_n, C_n)$, and so $|Z_{12}^{n}/\sqrt{Z_{11}^{n} Z_{22}^{n}}| \overset{D}{=} A_n/ \sqrt{B_n C_n}$. Statement (a) thus follows. Similarly, statement (b) follows from Cauchy-Schwarz inequality and from Proposition \ref{prop:identicaldistribution}.

The non-negativity of the terms in statements (c) and (d) comes from Proposition \ref{prop:identicaldistribution} and Jensen's inequality as is well known. For $n\geq 2$, $Z_{11}^{n}/n = 0$ implies  $X_1 = X_2 = \dots = X_n$, where $X_1, X_2, \dots, X_n$ are mutually independent standard Gaussian random variables as defined in Section \ref{sec:introduction}. We immediately note that the probability of the event $\{X_1 = X_2 = \dots = X_n\}$ is $0$. Thus, $Z_{11}^{n}/n > 0$ a.s., and hence $B_n > 0$ a.s.. Similarly, $C_n >0$ a.s. This proves statement (c). Finally, by Proposition \ref{prop:identicaldistribution}, $B=0$ implies that $W_{1}(t)$ is a constant on the interval $[0,1]$, the probability of which is $0$. Hence $B>0$ a.s. and similarly, $C>0$ a.s. This proves statement (d).
\end{proof}

Let $\phi_{B}(s) := E[e^{-sB}]$ and $\phi_{C}(s) := E[e^{-sC}]$. $\phi_{B}(s)$ and $\phi_{C}(s)$ are the moment generating functions of $B$ and $C$ respectively. Similarly, let $\phi_{B_n}(s)$ and $\phi_{C_n}(s)$ be the moment generating functions of $B_n$ and $C_n$ respectively. These functions can be computed explicitly, as the following lemma shows.
\begin{lemma} \label{lm:mgf}
We have
\begin{eqnarray*}
&&\phi_{B_n}(s) = \phi_{C_n}(s) =  \left( d_{n}\left(-2s/n\right) \right)^{-1/2}, \quad \text{for every $n \in \mathds{N}_{+}$}, \\
&& \phi_{B}(s) = \phi_{C}(s) = \left( \frac{\sinh\sqrt{2s}}{\sqrt{2s}} \right)^{-1/2}. 
\end{eqnarray*}
\end{lemma}
\begin{proof}
By Theorem \ref{thm1}, the joint moment generating function of $(Z_{11}^{n}, Z_{12}^{n}, Z_{22}^{n})$ is
\begin{eqnarray*}
&& E\left[\exp \left\{ -\frac{1}{2} (s_{11} Z_{11}^{n} + 2 s_{12} Z_{12}^{n} + s_{22} Z_{22}^{n}) \right\}  \right] 
= (d_n(\alpha(s_{11}, s_{12}, s_{22})) d_n(\beta(s_{11}, s_{12}, s_{22})))^{-1/2}. 
\end{eqnarray*}
Plugging $s_{11} = 2s/n$, $s_{12}=0$ and $s_{22}=0$ into the last display, it follows that
\begin{equation*}
E[e^{-s\,Z_{11}^{n} /n}] = \left( d_{n}(-2s/n) d_{n}(0)\right)^{-1/2}
 = \left( d_{n}(-2s/n) \right)^{-1/2}, 
\end{equation*}
where the last equality comes from the fact that $d_{n}(0)=1$. Note that since (by Proposition \ref{prop:identicaldistribution})  $Z_{11}^{n}/n \overset{D}{=} B_n$,
\begin{equation*}
\phi_{B_n}(s) = E[e^{-s B_n}] = E[e^{-s\,Z_{11}^{n} /n}] = \left( d_{n}(-2s/n) \right)^{-1/2}.  
\end{equation*}
Symmetrically, $\phi_{C_n}(s) =  \left( d_{n}\left(-2s/n\right) \right)^{-1/2}$.

Combining the results of Section 4.1 of \cite{ernst2019distribution} with Proposition \ref{prop:identicaldistribution}, the joint moment generating function of $(A,B,C)$ is
\begin{eqnarray*}
E\left[\exp \left\{ -\frac{1}{2} (s_{11} B + 2 s_{12} A + s_{22} C) \right\}  \right]
= \left( \frac{\sinh\sqrt{-\alpha} \sinh\sqrt{-\beta} }{\sqrt{-\alpha} \sqrt{-\beta}} \right)^{-1/2},  
\end{eqnarray*}
where $\alpha$ and $\beta$ are defined in \eqref{eq:alpha} and \eqref{eq:beta}. Plugging $s_{11} = 2s$, $s_{12}=0$ and $s_{22}=0$ into the last display, and recalling that $\sinh x/x$ equals $1$ at $x=0$, it follows that the moment generating function $\phi_{B}$ of $B$ is
\begin{equation*}
\phi_{B}(s) = E[e^{-sB}] = \left( \frac{\sinh\sqrt{2s}}{\sqrt{2s}} \right)^{-1/2}. 
\end{equation*}
Symmetrically, $\phi_{C}(s) = (\sinh\sqrt{2s} / \sqrt{2s})^{-1/2}$.
\end{proof}

We now proceed to give upper bounds for the second moments of $A_n -A$, $B_n -B$ and $C_n -C$.
\begin{proposition} \label{prop:secondmoment}
For $n>2$, we have $E[(A_n - A)^2] = \frac{5}{72} n^{-2} - \frac{7}{120} n^{-4}$. Hence, $E[(A_n - A)^2] \leq \frac{5}{72} n^{-2}$ for $n>2$.
\end{proposition}
\begin{proof}
From the definition \eqref{emenn}, a routine calculation shows that $M$ is a sublinear function in both its variables, with Lipshitz constant 1: 
$|M( s_2, t_2 ) - M( s_1, t_1 )| \leq \max(s_2 - s_1, t_2 - t_1)$, for $0\leq s_1 \leq s_2 \leq 1$ and $0 \leq t_1 \leq t_2 \leq 1$. It follows immediately that $$|M_{n}(s,t) - M(s,t) | \leq \frac{1}{n}.$$ By definitions of $A_n$ and $A$,
\begin{eqnarray*}
A_n - A &=& \int_{0}^{1} \int_{0}^{t} \left(M_{n}(s,t) - M(s,t) \right) \, dW_{1}(s) \, dW_{2}(t)  \\
&& + \int_{0}^{1} \int_{0}^{s} \left(M_{n}(s,t) - M(s,t) \right) \, dW_{2}(t) \, dW_{1}(s).  
\end{eqnarray*}
By Jensen's inequality,
\begin{eqnarray*}
(A_n - A)^{2} &\leq& 2 \left( \int_{0}^{1} \int_{0}^{t} \left(M_{n}(s,t) - M(s,t) \right) \, dW_{1}(s) \, dW_{2}(t) \right)^2 \\
&& + 2 \left( \int_{0}^{1} \int_{0}^{s} \left(M_{n}(s,t) - M(s,t) \right) \, dW_{2}(t) \, dW_{1}(s) \right)^2.  
\end{eqnarray*}
Taking expectations on both sides yields
\begin{eqnarray*}
&& E[(A_n - A)^{2}]  \\
&\leq& 2 E \left[ \left( \int_{0}^{1} \int_{0}^{t} \left(M_{n}(s,t) - M(s,t) \right) \, dW_{1}(s) \, dW_{2}(t) \right)^2 \right]   \\
&& + 2 E \left[ \left( \int_{0}^{1} \int_{0}^{s} \left(M_{n}(s,t) - M(s,t) \right) \, dW_{2}(t) \, dW_{1}(s) \right)^2 \right]  \\
&=& 2  \int_{0}^{1} \int_{0}^{t} E \left[ \left(M_{n}(s,t) - M(s,t) \right)^2 \right] \, ds\,  dt  
 + 2  \int_{0}^{1} \int_{0}^{s} E\left[ \left(M_{n}(s,t) - M(s,t) \right)^2  \right] \, dt\,  ds  \\
&\leq& 2  \int_{0}^{1} \int_{0}^{t} \frac{1}{n^2} \, ds\,  dt  
+ 2 \int_{0}^{1} \int_{0}^{s} \frac{1}{n^2} \, dt\,  ds  
= \frac{2}{n^2},  
\end{eqnarray*}
where in the first equality the It\^{o} isometry has been applied. In the last display, letting $n \rightarrow \infty$ yields $\lim_{n \rightarrow \infty} E[(A_n - A)^{2}] = 0$. Hence, $\lim_{n \rightarrow \infty} E[A_n^2] = E[A^2]$.

By \eqref{eq:An} in the appendix, we have
\begin{equation*}
A_n  =  \sum_{j,k=1}^{n} M\left( \frac{j-1}{n}, \frac{k-1}{n} \right) \left( W_1\left( \frac{j}{n}\right) -  W_1\left( \frac{j-1}{n}\right) \right) \left( W_2\left( \frac{k}{n}\right) -  W_2\left( \frac{k-1}{n}\right) \right).
\end{equation*}
Then
\begin{eqnarray}
&& E[A_n^2] \notag \\
&=&  E\Bigg[ \sum_{j,k,i,l=1}^{n} M\left( \frac{j-1}{n}, \frac{k-1}{n} \right) M\left( \frac{i-1}{n}, \frac{l-1}{n} \right) \left( W_{1}\left(\frac{j}{n}\right) - W_{1}\left(  \frac{j-1}{n}\right) \right) \notag \\
&& \left( W_{2}\left(\frac{k}{n}\right) - W_{2}\left(  \frac{k-1}{n}\right) \right)  
\left( W_{1}\left(\frac{i}{n}\right) - W_{1}\left(  \frac{i-1}{n}\right) \right)
\left( W_{2}\left(\frac{l}{n}\right) - W_{2}\left(  \frac{l-1}{n}\right) \right)  \Bigg]  \notag \\
&=& \sum_{j,k,i,l=1}^{n} \Bigg\{ M\left( \frac{j-1}{n}, \frac{k-1}{n} \right) M\left( \frac{i-1}{n}, \frac{l-1}{n} \right) \notag \\
&&   \quad \quad \quad \quad \times E\left[ \left( W_{1}\left(\frac{j}{n}\right) - W_{1}\left(  \frac{j-1}{n}\right) \right) \left( W_{1}\left(\frac{i}{n}\right) - W_{1}\left(  \frac{i-1}{n}\right) \right) \right]  \notag \\
&&   \quad \quad \quad \quad  \times E\left[ \left( W_{2}\left(\frac{k}{n}\right) - W_{2}\left(  \frac{k-1}{n}\right) \right) \left( W_{2}\left(\frac{l}{n}\right) - W_{2}\left(  \frac{l-1}{n}\right) \right)  \right]  \Bigg\} \notag \\
&=& \sum_{j,k,i,l=1}^{n}  M\left( \frac{j-1}{n}, \frac{k-1}{n} \right) M\left( \frac{i-1}{n}, \frac{l-1}{n} \right) \cdot \frac{1}{n} \mathds{1}_{\{i=j\}} \cdot \frac{1}{n} \mathds{1}_{\{l=k\}} \notag \\
&=& \sum_{j,k=1}^{n} M\left( \frac{j-1}{n}, \frac{k-1}{n} \right)^2 \cdot \frac{1}{n^2} \notag \\
&=& 2 \sum_{j=1}^{n}\sum_{k=1}^{j-1} \frac{1}{n^2} \left( \frac{k-1}{n} - \frac{(j-1)(k-1)}{n^2} \right)^2 
   + \sum_{j=1}^{n} \frac{1}{n^2} \left( \frac{j-1}{n} - \frac{(j-1)^2}{n^2} \right)^2,  \label{eq:longexpandofAn}
\end{eqnarray}
where the second equality follows by the independence of the Wiener processes $W_{1}$ and $W_{2}$. Note that the first term on the right-hand side of \eqref{eq:longexpandofAn} is $2$ times a double summation of a polynomial of $j$ and $k$, calculating the summation with respect to $k$ by Faulhaber's formula yields $2$ times a single summation of a polynomial of $j$ over $j=1,\dots, n$. Applying Faulhaber's formula again to this summation, we have that the first term on the right-hand side of \eqref{eq:longexpandofAn} is equal to
\begin{equation}
\frac{1}{90} - \frac{1}{30} n^{-1} + \frac{1}{36} n^{-2} - \frac{7}{180} n^{-4} + \frac{1}{30} n^{-5}.  \label{eq:Andoublesum}
\end{equation}
 Note that the second term on the right-hand side of \eqref{eq:longexpandofAn} is a single summation of a polynomial of $j$, applying Faulhaber's formula again gives
\begin{equation}
 \frac{1}{30} n^{-1} - \frac{1}{30} n^{-5}.  \label{eq:Ansinglesum}
\end{equation}
 Combining \eqref{eq:longexpandofAn}, \eqref{eq:Andoublesum} and \eqref{eq:Ansinglesum}, we have
\begin{equation}
 E[A_{n}^{2}] =  \frac{1}{90} + \frac{1}{36} n^{-2}  - \frac{7}{180} n^{-4}.  \label{eq:expectationofAn2}
\end{equation}
In the last display, letting $n \rightarrow \infty$, we have
\begin{equation}
E[A^2] = \lim_{n \rightarrow \infty} E[A_{n}^{2}] =  \frac{1}{90}. \label{eq:expectationofA2}
\end{equation}

In what follows, we proceed to calculate the expectation of $A_n A$, which, of course, is handy to compute the variance of $A_n-A$. By Fubini's theorem and the independence of $W_1$ and $W_2$,
\begin{eqnarray*}
&& E\left[ \left( W_{1}\left(\frac{j}{n}\right) - W_{1}\left(  \frac{j-1}{n}\right) \right)
\left( W_{2}\left(\frac{k}{n}\right) - W_{2}\left(  \frac{k-1}{n}\right) \right)  \int_{0}^{1} W_{1}(t) W_{2}(t) \, dt \right]  \\
&=& \int_{0}^{1} E\left[ \left( W_{1}\left(\frac{j}{n}\right) - W_{1}\left(  \frac{j-1}{n}\right) \right)
\left( W_{2}\left(\frac{k}{n}\right) - W_{2}\left(  \frac{k-1}{n}\right) \right) W_{1}(t) W_{2}(t) \right] \, dt  \\
&=& \int_{0}^{1} E\left[ \left( W_{1}\left(\frac{j}{n}\right) - W_{1}\left(  \frac{j-1}{n}\right) \right) W_{1}(t) \right]
 E\left[ \left( W_{2}\left(\frac{k}{n}\right) - W_{2}\left(  \frac{k-1}{n}\right) \right) W_{2}(t) \right] \, dt  \\
&=& \int_{0}^{1} \left( t \wedge \left( \frac{j}{n} \right)  - t \wedge \left( \frac{j-1}{n} \right)\right) \cdot 
\left( t \wedge \left( \frac{k}{n} \right)  - t \wedge \left( \frac{k-1}{n} \right)\right) \, dt  \\
&=&  \frac{1}{n^2} + \frac{1}{2} \frac{1}{n^3} - \frac{j \vee k}{n^3} - \frac{1}{6} \frac{1}{n^3} \mathds{1}_{\{j=k\}}.
\end{eqnarray*} 
Similarly,
\begin{eqnarray*}
&& E \left[ \left( W_{1}\left(\frac{j}{n}\right) - W_{1}\left(  \frac{j-1}{n}\right) \right)
\left( W_{2}\left(\frac{k}{n}\right) - W_{2}\left(  \frac{k-1}{n}\right) \right)  \int_{0}^{1} W_{1}(t) \, dt \int_{0}^{1} W_{1}(t)  \, dt \right]  \\
&=& \int_{0}^{1} \int_{0}^{1} E\left[ \left( W_{1}\left(\frac{j}{n}\right) - W_{1}\left(  \frac{j-1}{n}\right) \right)
\left( W_{2}\left(\frac{k}{n}\right) - W_{2}\left(  \frac{k-1}{n}\right) \right) W_{1}(t) W_{2}(s) \right] \,ds\,dt \\
&=& \int_{0}^{1} \int_{0}^{1} E\left[ \left( W_{1}\left(\frac{j}{n}\right) - W_{1}\left(  \frac{j-1}{n}\right) \right) W_{1}(t) \right]
 E\left[ \left( W_{2}\left(\frac{k}{n}\right) - W_{2}\left(  \frac{k-1}{n}\right) \right) W_{2}(s) \right] \,ds\,dt \\
&=& \int_{0}^{1} \int_{0}^{1} \left( t \wedge \left( \frac{j}{n} \right)  - t \wedge \left( \frac{j-1}{n} \right)\right) \cdot 
\left( s \wedge \left( \frac{k}{n} \right)  - s \wedge \left( \frac{k-1}{n} \right)\right) \,ds\,dt \\
&=& \left( \frac{1}{n} + \frac{1}{2} \frac{1}{n^2} - \frac{j}{n^2} \right) \cdot \left( \frac{1}{n} + \frac{1}{2} \frac{1}{n^2} - \frac{k}{n^2} \right).
\end{eqnarray*}
Combing the last two displays and by linearity of expectation, we have
\begin{eqnarray}
&& E[A_n A] \notag \\
&=& \sum_{j,k=1}^{n} M\left( \frac{j-1}{n}, \frac{k-1}{n} \right) \bigg[ \frac{1}{n^2} + \frac{1}{2} \frac{1}{n^3} - \frac{j \vee k}{n^3} - \frac{1}{6} \frac{1}{n^3} \mathds{1}_{\{j=k\}}   \notag \\
&& \quad \quad \quad - \left( \frac{1}{n} + \frac{1}{2} \frac{1}{n^2} - \frac{j}{n^2} \right) \cdot \left( \frac{1}{n} + \frac{1}{2} \frac{1}{n^2} - \frac{k}{n^2} \right) \bigg]  \notag \\
&=&  2 \sum_{j=1}^{n} \sum_{k=1}^{j-1} \left( \frac{k-1}{n} - \frac{j-1}{n}\cdot \frac{k-1}{n}\right) \bigg[ \frac{1}{n^2} + \frac{1}{2} \frac{1}{n^3} - \frac{j }{n^3}  \notag \\
&& \quad \quad \quad \quad \quad - \left( \frac{1}{n} + \frac{1}{2} \frac{1}{n^2} - \frac{j}{n^2} \right) \cdot \left( \frac{1}{n} + \frac{1}{2} \frac{1}{n^2} - \frac{k}{n^2} \right) \bigg]  \notag \\
&& + \sum_{j=1}^{n} \left( \frac{j-1}{n} - \frac{j-1}{n}\cdot \frac{j-1}{n}\right) \bigg[ \frac{1}{n^2} + \frac{1}{2} \frac{1}{n^3} - \frac{j }{n^3}  - \frac{1}{6} \frac{1}{n^3} \notag \\
&& \quad \quad \quad  - \left( \frac{1}{n} + \frac{1}{2} \frac{1}{n^2} - \frac{j}{n^2} \right) \cdot \left( \frac{1}{n} + \frac{1}{2} \frac{1}{n^2} - \frac{j}{n^2} \right) \bigg].  \label{eq:longexpandofAnA}
\end{eqnarray}
By a similar argument to that of the calculation of the right-hand side of \eqref{eq:longexpandofAn}, the right-hand side of \eqref{eq:longexpandofAnA} is
\begin{equation*}
\frac{1}{90} - \frac{1}{48}n^{-2} + \frac{7}{720}n^{-4}.
\end{equation*}
Hence,
\begin{equation*}
E[A_n A] = \frac{1}{90} - \frac{1}{48}n^{-2} + \frac{7}{720}n^{-4}.
\end{equation*}
Together with \eqref{eq:expectationofAn2} and \eqref{eq:expectationofA2}, we have
\begin{equation*}
E[(A_n - A)^2] = E[A_{n}^{2}] + E[A^2] - 2 E[A_n A] = \frac{5}{72} n^{-2} - \frac{7}{120} n^{-4}.
\end{equation*}
\end{proof}

\begin{proposition}  \label{prop:upperboundforBn-B}
We have $E[(B_n- B)^2] = E[(C_n - C)^2]= \frac{5}{36}n^{-2} - \frac{4}{45} n^{-4}$. Hence, $E[(B_n- B)^2] = E[(C_n - C)^2] \leq \frac{5}{36}n^{-2}$.
\end{proposition}
\begin{proof}
The second assertion is a direct result of the first, which we proceed to establish.
By the definitions of $B_n$ and $B$,
\begin{eqnarray}
B_n - B &=& 2\int_{0}^{1} \int_{0}^{t} \left(M_{n}(s,t) - M(s,t) \right) \, dW_{1}(s) \, dW_{1}(t)  \notag \\
&& + \frac{1}{n} \sum_{j=1}^{n} M\left(\frac{j-1}{n}, \frac{j-1}{n}\right) - \int_{0}^{1} M(t,t)\, dt. \label{eq:BnminusB}
\end{eqnarray}
By a standard property of the Wiener stochastic integral,
\begin{equation*}
E \left[ \int_{0}^{1} \int_{0}^{t} \left(M_{n}(s,t) - M(s,t) \right) \, dW_{1}(s) \, dW_{1}(t) \right] =0. 
\end{equation*}
Taking squares and then expectation on both sides of \eqref{eq:BnminusB}, we have (after rearrangement of terms) that
\begin{eqnarray}
E[ (B_n -B)^2 ] &=& 4 E\left[ \left(\int_{0}^{1} \int_{0}^{t} \left(M_{n}(s,t) - M(s,t) \right) \, dW_{1}(s) \, dW_{1}(t)\right)^2 
\right] \notag \\
&& +  \left( \frac{1}{n} \sum_{j=1}^{n} M\left(\frac{j-1}{n}, \frac{j-1}{n}\right) - \int_{0}^{1} M(t,t)\, dt \right)^2. \label{eq:boundB2} 
\end{eqnarray}
By the It\^{o} isometry, the first term on the right-hand side of \eqref{eq:boundB2} is
\begin{eqnarray}
&& 4 \int_{0}^{1} \int_{0}^{t} E \left[ \left(M_{n}(s,t) - M(s,t) \right)^2 \right] \, ds\,  dt  \notag \\
&=& 4 \int_{0}^{1} \int_{0}^{t}  \left(M_{n}(s,t) - M(s,t) \right)^2  \, ds\,  dt  
= 4 \sum_{j=1}^{n} \int_{\frac{j-1}{n}}^{\frac{j}{n}} \int_{0}^{t}  \left(M_{n}(s,t) - M(s,t) \right)^2  \, ds\,  dt \notag  \\
&=& 4 \sum_{j=1}^{n} \int_{\frac{j-1}{n}}^{\frac{j}{n}} \left( \sum_{k=1}^{j-1} \int_{\frac{k-1}{n}}^{\frac{k}{n}}  \left(M_{n}(s,t) - M(s,t) \right)^2  \, ds  + \int_{\frac{j-1}{n}}^{t} \left(M_{n}(s,t) - M(s,t) \right)^2  \, ds \right) \,dt  \notag \\
&=& 4 \sum_{j=1}^{n} \sum_{k=1}^{j-1} \int_{\frac{j-1}{n}}^{\frac{j}{n}} \int_{\frac{k-1}{n}}^{\frac{k}{n}} \left(M_{n}(s,t) - M(s,t) \right)^2 \,ds\,dt \notag \\
 &&  + 4 \sum_{j=1}^{n} \int_{\frac{j-1}{n}}^{\frac{j}{n}} \int_{\frac{j-1}{n}}^{t} \left(M_{n}(s,t) - M(s,t) \right)^2 \,ds\,dt  \notag \\
&=& 4 \sum_{j=1}^{n} \sum_{k=1}^{j-1} \int_{\frac{j-1}{n}}^{\frac{j}{n}} \int_{\frac{k-1}{n}}^{\frac{k}{n}} \left(\frac{k-1}{n}- \frac{k-1}{n}\cdot \frac{j-1}{n} -s + st \right)^2 \,ds\,dt \notag \\
&& + 4 \sum_{j=1}^{n} \int_{\frac{j-1}{n}}^{\frac{j}{n}} \int_{\frac{j-1}{n}}^{t} \left( \frac{j-1}{n} - \frac{j-1}{n} \cdot \frac{j-1}{n} -s+st \right)^2 \, ds \, dt  \notag \\
&=&  4 \sum_{j=1}^{n} \sum_{k=1}^{j-1} \left( \frac{2j^2 + 2k^2 + 3jk}{6n^6} - \frac{(4n+5)j}{6n^6} - \frac{(3n+5)k}{6n^6} + \frac{6n^2 + 15n +11}{18n^6} \right)  \notag \\
&& + 4 \sum_{j=1}^{n} \left( \frac{7 j^2}{12n^6} - \frac{5(n+2)j}{12n^6} + \frac{15n^2 + 51n +55}{180 n^6} \right), \label{eq:longexpansionofBn1}
\end{eqnarray}
where in the last equality we have explicitly calculated the two double integrals. Note that the first term on the right-hand side of \eqref{eq:longexpansionofBn1} is $4$ times a double summation of a polynomial of $j$ and $k$, calculating the summation with respect to $k$ by Faulhaber's formula yields $4$ times a single summation of a polynomial of $j$ over $j=1,\dots, n$. Applying Faulhaber's formula again to this summation, we have that the first term on the right-hand side of \eqref{eq:longexpansionofBn1} is
\begin{equation}
\frac{5}{36} n^{-2} - \frac{5}{18} n^{-3} + \frac{1}{12}n^{-4} + \frac{1}{18} n^{-5}. \label{eq:momentofBdoublesum}
\end{equation}
Again, by Faulhaber's formula, the second term on the right-hand side of \eqref{eq:longexpansionofBn1} is
\begin{equation}
\frac{5}{18}n^{-3} - \frac{1}{5} n^{-4} - \frac{1}{18} n^{-5}. \label{eq:momentofBsinglesum}{}
\end{equation}
Combing \eqref{eq:momentofBdoublesum} and \eqref{eq:momentofBsinglesum}, the first term on the right-hand side of \eqref{eq:boundB2} is
\begin{equation}
\frac{5}{36}n^{-2} - \frac{7}{60} n^{-4}. \label{eq:B2term1}
\end{equation}
The second term on the right-hand side of \eqref{eq:boundB2} is
\begin{eqnarray}
&& \left( \frac{1}{n} \sum_{j=1}^{n} \left( \frac{j-1}{n} - \left(\frac{j-1}{n}\right)^2 \right) - \int_{0}^{1} (t - t^2)\, dt \right)^2 \notag  \\
&=&  \left( \frac{n^2-1}{6n^2} - \frac{1}{6} \right)^2 = \frac{1}{36} n^{-4}. \label{eq:B2term2}  
\end{eqnarray}
Combining \eqref{eq:boundB2}, \eqref{eq:B2term1}, \eqref{eq:B2term2} gives $E[(B_n- B)^2] = \frac{5}{36}n^{-2} - \frac{4}{45} n^{-4}$. Symmetrically, $E[(C_n- C)^2] = \frac{5}{36}n^{-2} - \frac{4}{45} n^{-4}$ too. This completes the proof.
\end{proof}

\subsection{An upper bound for the Wasserstein distance.}
\label{sec:UBWass}
In this section, we will derive an upper bound for $d_{W}(A_{n}/ \sqrt{B_{n} C_{n}}, A/\sqrt{BC})$. The result relies on three preparatory lemmas. The first, Lemma \ref{lm:momentinverse}, is a special case of Proposition 1 in \cite{ernst2019distribution}. When used in conjunction with Lemma \ref{lm:mgf}, it shows that we must have a good lower-bound handle on the behavior of $d_n$, which is the topic of the Lemma \ref{lm:lowerboundfordn}. These then culminate in showing (Lemma \ref{lm:boundforBinverse}) that $B$ and $B_n$ have inverse moments, with the latter being uniformly bounded in $n$. This fact may seem surprising, since, as second chaos variables, negative moments can explode, but this does not apply because $B,B_n$ are non-centered, and a.s. positive. The uniformity over $n$ in Lemma \ref{lm:boundforBinverse} is a consequence of the convergence of the moment-generating functions of the $B_n$'s to a limit which decays rapidly at $+\infty$ (at the rate $\sqrt{2s}e^{-2s}$), ensuring control of the tails. 

We exploit the explicit nature of these expressions to prove Lemma \ref{lm:boundforBinverse} and the results that precede it, but our strategy could also work for other processes, for instance by invoking dominated convergence and by controlling $d_n$ via its constituent eigenvalues. This means that our methodology could handle other processes, or other quadratic forms than $B_n$, if one could still control $d_n$, via the properties of the matrix $K_n$, whose positive-definite character is very general. This is an important point in understanding the ingredients in the proof of Lemma \ref{lm:lowerboundfordn}. We obtain lower bounds for $d_n$ by estimating selected terms in its sum representation, ignoring others because none of them are negative, and the positive-definite property of $K_n$ is the reason all terms in the sum are non-negative. This last justification is not entirely trivial, and though it is not used in our proofs because all our formulas are explicit, it is worth mentioning the reason here which is generically true. We are interested in lower bounds on the moment-generating function of $B_n$, which equals $d_n(-2s/n)=\prod_{k=2}^n(1+2\lambda_k s/n)$. Since all $\lambda_k$ are positive, this expression is thus a polynomial in $s$ with positive coefficients. That positivity translates into the one used in the proof of Lemma \ref{lm:lowerboundfordn}. 
\begin{lemma} \label{lm:momentinverse}
Let $X$ be a random variable satisfying $X>0$ a.s. and $\phi_{X}(s) = E[e^{-sX}]$ be its moment generating function. Then for every $m \in \mathds{N}_{+}$,
\begin{equation*}
E\left[ X^{-m} \right] = \frac{1}{(m-1)!} \int_{0}^{\infty} s^{m-1} \phi_{X}(s) \, ds. 
\end{equation*}
\end{lemma}
\begin{proof}
See Proposition 1 in \cite{ernst2019distribution}.
\end{proof}
We now turn to Lemma \ref{lm:lowerboundfordn}.
\begin{lemma}  \label{lm:lowerboundfordn}
For $n \geq 11$, we have
\begin{enumerate}[(a)]
\item $d_{n}(-2s/n) \geq 1$ for $s \geq 0$;
\item $d_{n}(-2s/n) \geq 2^{5} {n \choose 11}\, n^{-11} s^{5} $ for $s \geq 0$;
\item $d_{n}(-2s/n) \geq \left( e^{\sqrt{s/2} } - e^{-\sqrt{s/2} } \right)/ \sqrt{10s}$ for $ 0 \leq s \leq n^2 /2$.
\end{enumerate}
\end{lemma}
\begin{proof}
It follows from \eqref{eq:representationfordinmaindocument} that
\begin{eqnarray}
d_{n}\left( - \frac{2s}{n} \right) &=& \frac{(-1)^{n-1}}{n \cdot 2^{n-1}} \, \sum_{k=1}^{\lceil n/2 \rceil} {n \choose 2k-1} \left(-\frac{2s}{n^2} -2\right)^{n-(2k-1)} \left( \left(-\frac{2s}{n^2} -2\right)^2 -4 \right)^{k-1} \notag  \\
&=&  \sum_{k=1}^{\lceil n/2 \rceil} \frac{1}{n} \, {n \choose 2k-1} \left( \frac{s}{n^2} +1 \right)^{n-(2k-1)}  \left( \frac{2s}{n^2} + \frac{s^2}{n^4} \right)^{k-1}. \label{eq:partialsum} 
\end{eqnarray}
Note the first term of the summation on the right-hand side of \eqref{eq:partialsum} is $(s/n^2 +1)^{n-1}$. Then,
\begin{equation*}
d_{n}\left( - \frac{2s}{n} \right) \geq \left( \frac{s}{n^2} +1 \right)^{n-1} \geq 1. 
\end{equation*}
This proves statement (a). We now note that $n \geq 11$, $\lceil n/2 \rceil \geq 6$. Let us consider the sixth term of the summation on the right-hand side of \eqref{eq:partialsum}, i.e.
\begin{equation*}
\frac{1}{n} \, {n \choose 11} \left( \frac{s}{n^2}  +1\right)^{n-11} \left( \frac{2s}{n^2} + \frac{s^2}{n^4} \right)^{5}
\geq \frac{1}{n} \, {n \choose 11} \left( \frac{2s}{n^2} \right)^{5}
=  2^{5} {n \choose 11}\, n^{-11} s^{5}. 
\end{equation*}
Then $d_{n}(-2s/n) \geq 2^{5} {n \choose 11}\, n^{-11} s^{5} $. This completes the proof of statement (b). Finally, it follows from \eqref{eq:representationfordinmaindocument} that
\begin{eqnarray}
d_{n}\left( - \frac{2s}{n} \right) = 
\frac{1}{2 \sqrt{2s + s^2 /n^2}} \left[ \left( \frac{s}{n^2} + 1 + \sqrt{\frac{2s}{n^2}  + \frac{s^2}{n^4}} \right)^{n} - \left( \frac{s}{n^2} + 1 - \sqrt{\frac{2s}{n^2}  + \frac{s^2}{n^4}} \right)^{n} \right]. \label{eq:assertionc1} 
\end{eqnarray}
Recall the useful fact that $\log(1+x) \geq x/2$ for $0 \leq x \leq 1$. Then for $ 0 \leq s \leq n^2 /2$,
\begin{eqnarray}
\left( \frac{s}{n^2} + 1 + \sqrt{\frac{2s}{n^2}  + \frac{s^2}{n^4}} \right)^{n}
\geq  \left( 1+ \frac{\sqrt{2s}}{n} \right)^{n}
= e^{n \log(1 + \sqrt{2s}/n)}
\geq e^{\sqrt{s/2}}, \label{eq:assertionc2} 
\end{eqnarray}
Further,
\begin{eqnarray}
\left( \frac{s}{n^2} + 1 - \sqrt{\frac{2s}{n^2}  + \frac{s^2}{n^4}} \right)^{n}
= \left( \frac{s}{n^2} + 1 + \sqrt{\frac{2s}{n^2}  + \frac{s^2}{n^4}} \right)^{- n}
\leq e^{-\sqrt{s/2}}. \label{eq:assertionc3} 
\end{eqnarray}
Note that
\begin{eqnarray*}
2 s + \frac{s^2}{n^2} \leq 2 s + \frac{n^2}{2} \cdot \frac{s}{n^2} = \frac{5s}{2}, 
\end{eqnarray*}
Together with \eqref{eq:assertionc1}, \eqref{eq:assertionc2} and \eqref{eq:assertionc3}, statement (c) follows.
\end{proof}
We now turn to Lemma \ref{lm:boundforBinverse} below.
\begin{lemma} \label{lm:boundforBinverse}
We have 
\begin{enumerate}[(a)]
\item $E[B^{-m}] = E[C^{-m}] < \infty$ for every $m \in \mathds{N}_{+}$;
\item $\sup_{n\geq 11}  E[B_{n}^{-1}] = \sup_{n\geq 11}  E[C_{n}^{-1}] < \infty$.
\end{enumerate}
\end{lemma}
\begin{proof}
We first consider statement (a). Since $B$ and $C$ are identically distributed, it suffices to prove that $E[B^{-m}]  < \infty$ for every $m \in \mathds{N}_{+}$. Applying Lemma \ref{lm:momentinverse} gives
\begin{equation*}
E[B^{-m}] = \frac{1}{(m-1)!} \int_{0}^{\infty} s^{m-1} \phi_{B}(s) \, ds. 
\end{equation*}
By Lemma \ref{lm:mgf}, we have $\phi_{B}(0)=1$ and $\phi_{B}(s) \sim 2^{3/2}\,  s^{1/4} \,  e^{-\sqrt{s/2}}$ as $s \rightarrow \infty$. Then the boundedness of $E[B^{-m}]$ follows immediately. For statement (b), similarly, we need only prove that $\sup_{n\geq 11}  E[B_{n}^{-1}]  < \infty$. By Lemma \ref{lm:mgf}, Lemma \ref{lm:momentinverse} and Lemma \ref{lm:lowerboundfordn}, we have
\begin{eqnarray*}
&&E[B_{n}^{-1}] = \int_{0}^{\infty}  \phi_{B_n}(s) \, ds
= \int_{0}^{\infty}   d_{n}\left(- \frac{2s}{n}\right)^{-1/2} \, ds \\
&=& \int_{0}^{1}   d_{n}\left(- \frac{2s}{n}\right)^{-1/2} \, ds
   + \int_{1}^{n^2/2}   d_{n}\left(- \frac{2s}{n}\right)^{-1/2} \, ds
   + \int_{n^2/2}^{\infty}   d_{n}\left(- \frac{2s}{n}\right)^{-1/2} \, ds \\
&\leq& \int_{0}^{1}  1 \, ds
 + \int_{1}^{n^2/2}   \left( \left( e^{\sqrt{s/2} } - e^{-\sqrt{s/2} } \right)/ \sqrt{10s} \right)^{-1/2} \, ds  
  +  \int_{n^2/2}^{\infty}  2^{-\frac{5}{2}} {n \choose 11}^{-\frac{1}{2}} n^{\frac{11}{2}} s^{-\frac{5}{2}} \, ds  \\
&=& 1 + \int_{1}^{n^2/2}   \left( \left( e^{\sqrt{s/2} } - e^{-\sqrt{s/2} } \right)/ \sqrt{10s} \right)^{-1/2} \, ds 
  + \frac{1}{3} {n \choose 11}^{-\frac{1}{2}} n^{\frac{5}{2}},  
\end{eqnarray*}
where in the first inequality we have applied Lemma \ref{lm:lowerboundfordn}. Taking the supremum over $n\geq 11$ on the both sides of the last display yields
\begin{equation}
\sup_{n\geq 11}  E[B_{n}^{-1}]
 \leq 1 + \int_{1}^{\infty}   \left( \left( e^{\sqrt{s/2} } - e^{-\sqrt{s/2} } \right)/ \sqrt{10s} \right)^{-1/2} \, ds
  + \sup_{n\geq 11} \frac{1}{3} {n \choose 11}^{-\frac{1}{2}} n^{\frac{5}{2}}.  \label{eq:boundednessofbniverse}  
\end{equation}
The boundedness of the second term on the right-hand side of \eqref{eq:boundednessofbniverse} follows by the fact that, as $n \rightarrow \infty$, $\left( \left( e^{\sqrt{s/2} } - e^{-\sqrt{s/2} } \right)/ \sqrt{10s} \right)^{-1/2} \sim (10s)^{1/4}\, e^{-\sqrt{s/8}}$. The boundedness of the third term on the right-hand side of \eqref{eq:boundednessofbniverse} follows by the fact ${n \choose 11}^{-\frac{1}{2}} n^{\frac{5}{2}} \rightarrow 0$ as $n \rightarrow \infty$. This completes the proof of statement (b).
\end{proof}

Let us define $C_{m} := E[B^{-m}] = E[C^{-m}]$ for $m=1,2,3$ and $C_4 := \sup_{n\geq 11}  E[B_{n}^{-1}] = \sup_{n\geq 11}  E[C_{n}^{-1}]$. With above preparations in hand, we are ready to reveal the section's main result.
\begin{theorem} \label{thm:boundwassertain}
For $n \geq 11$, with $A,B,C,A_n,B_n,C_n$ defined in Section \ref{sec:propertyofAnBnCnABC}, we have
\begin{equation*}
d_{W}(\theta_n, \theta) \leq E \left| \frac{A_n}{\sqrt{B_n C_n}} -  \frac{A}{\sqrt{B C}} \right| \leq \frac{C_5}{n},  
\end{equation*}
where, via the constants $C_1,C_3,C_4$ defined above, the constant $C_5$ is defined as
\begin{equation*}
C_5 := \frac{1}{12} \left\{ \frac{1}{132} \, \left[ \frac{5}{2} (C_3 +C_4) \right]^{\frac{1}{2}} + 2 \right\} 
\left[ \frac{5}{2} (C_3 +C_4) \right]^{\frac{1}{2}} + \frac{1}{6} \left(\frac{5}{2}\right)^{\frac{1}{2}} C_1.
\end{equation*}
\end{theorem}
\begin{proof}
By Proposition \ref{prop:identicaldistribution}, we have
\begin{equation*}
d_{W}(\theta_n, \theta) 
= d_{W} \left( \frac{A_n}{\sqrt{B_n C_n}}, \frac{A}{\sqrt{B C}}  \right)
= \sup_{f \in \mathrm{Lip}(1)} \left|Ef\left( \frac{A_n}{\sqrt{B_n C_n}} \right) - E f\left(\frac{A}{\sqrt{B C}}\right) \right|.
\end{equation*}
For every $f \in \mathrm{Lip}(1)$, and every pair of integrable random variables $(X,Y)$ on the same probability space,
\begin{equation*}
\left|Ef\left( X \right) - E f\left(Y\right) \right|
\leq E \left|  f\left( X \right) -  f\left( Y \right)\right| \\
\leq E \left| X - Y \right|. 
\end{equation*}
Taking the supremun over $f \in \mathrm{Lip}(1)$ on both sides of the above displays with $X,Y$ replaced by $A_n /\sqrt{B_n C_n}, A/\sqrt{B C}$ yields
\begin{equation*}
d_{W}(\theta_n, \theta) \leq E \left| \frac{A_n}{\sqrt{B_n C_n}} -  \frac{A}{\sqrt{B C}} \right|.  
\end{equation*}
Thus, we need only bound the expectation of $| A_{n}/ \sqrt{B_n C_n} - A/\sqrt{BC} |$. Note that
\begin{eqnarray*}
&&\frac{A_n}{\sqrt{B_n C_n}} -  \frac{A}{\sqrt{B C}}
= \frac{A_n \sqrt{B C} - A \sqrt{B_n C_n}}{\sqrt{B_n C_n} \sqrt{B C}}  \\
&=& \frac{A_n(\sqrt{B C} - \sqrt{B_n C_n}) + (A_n -A) \sqrt{B_n C_n}}{\sqrt{B_n C_n} \sqrt{B C}} \\
&=& \frac{A_n}{\sqrt{B_n C_n}} \cdot \frac{\sqrt{B C} - \sqrt{B_n C_n}}{\sqrt{B C}} + \frac{A_n -A}{\sqrt{B C}}. 
\end{eqnarray*}
Then,
\begin{eqnarray}
&&\left| \frac{A_n}{\sqrt{B_n C_n}} -  \frac{A}{\sqrt{B C}} \right|
\leq \left| \frac{A_n}{\sqrt{B_n C_n}}\right| \cdot \frac{|\sqrt{B C} - \sqrt{B_n C_n}|}{\sqrt{B C}}
+ \frac{|A_n -A|}{\sqrt{B C}}   \notag \\
&\leq& \frac{|\sqrt{B C} - \sqrt{B_n C_n}|}{\sqrt{B C}}
+ \frac{|A_n -A|}{\sqrt{B C}}  \notag   \\
&=& \frac{|(\sqrt{B}- \sqrt{B_n})\sqrt{C} + \sqrt{B_n}(\sqrt{C}- \sqrt{C_n})|}{\sqrt{B C}}
+ \frac{|A_n -A|}{\sqrt{B C}} \notag \\
&\leq& \frac{|\sqrt{B}- \sqrt{B_n}|}{\sqrt{B}} + \frac{\sqrt{B_n}}{\sqrt{B}}\, \frac{|\sqrt{C}- \sqrt{C_n}|}{\sqrt{C}} 
+ \frac{|A_n -A|}{\sqrt{B C}},  \label{eq:differenceoffraction} 
\end{eqnarray}
where in the second inequality we have invoked Corollary \ref{coro:1}. By the inequality of arithmetic and geometric means, 
\begin{eqnarray*}
\frac{|\sqrt{B}- \sqrt{B_n}|}{\sqrt{B}}
= \frac{|B_n - B|}{\sqrt{B} (\sqrt{B} + \sqrt{B_n})}
\leq \frac{|B_n - B|}{2 \sqrt{B} \sqrt{\sqrt{B}  \sqrt{B_n}}}
= \frac{|B_n - B|}{2 B^{\frac{3}{4}}\, B_{n}^{\frac{1}{4}}}. 
\end{eqnarray*}
Then by H\"{o}lder's inequality,
\begin{eqnarray}
&&E\left[ \frac{|\sqrt{B}- \sqrt{B_n}|}{\sqrt{B}} \right]
\leq E\left[ \frac{|B_n - B|}{2 B^{\frac{3}{4}}\, B_{n}^{\frac{1}{4}}} \right]  
\leq \left \{E\left[ \left(B_n - B\right)^2 \right] \right \} ^{\frac{1}{2}}
 \left \{ E\left[ \frac{1}{4} B^{-\frac{3}{2}} \, B_{n}^{-\frac{1}{2}} \right] \right \}^{\frac{1}{2}} \notag  \\
&\leq& \left\{ E\left[ \left(B_n - B\right)^2 \right] \right\}^{\frac{1}{2}}
 \left\{ E\left[ \frac{1}{8} \left( B^{-3} + B_{n}^{-1}\right) \right] \right\}^{\frac{1}{2}}  
\leq \left(\frac{5}{36n^2}\right)^{\frac{1}{2}} \, \left[ \frac{1}{8} (C_3 +C_4) \right]^{\frac{1}{2}}  \notag \\
&=& \frac{1}{12} \left[ \frac{5}{2} (C_3 +C_4) \right]^{\frac{1}{2}} \, \frac{1}{n}, \label{eq:difference1} 
\end{eqnarray}
where the third inequality follows from the inequality of arithmetic and geometric means, and the fourth inequality is due to Proposition \ref{prop:upperboundforBn-B}. Similarly, we have
\begin{equation*}
E\left[ \frac{|\sqrt{C}- \sqrt{C_n}|}{\sqrt{C}} \right] \leq \frac{1}{12} \left[ \frac{5}{2} (C_3 +C_4) \right]^{\frac{1}{2}} \, \frac{1}{n},
\end{equation*}
and
\begin{eqnarray*}
&& E\left[ \frac{\sqrt{B_n}}{\sqrt{B}} \right]
= E\left[ \frac{\sqrt{B_n} - \sqrt{B}}{\sqrt{B}} +1 \right]
\leq E\left[ \frac{|\sqrt{B_n} - \sqrt{B}|}{\sqrt{B}} +1 \right]  \\
&\leq& \frac{1}{12}  \left[ \frac{5}{2} (C_3 +C_4) \right]^{\frac{1}{2}} \, \frac{1}{n} +1 
\leq \frac{1}{132} \, \left[ \frac{5}{2} (C_3 +C_4) \right]^{\frac{1}{2}} + 1.  
\end{eqnarray*}
Since $(B_n, B)$ and $(C_n ,C)$ are independent, we have
\begin{eqnarray}
&& E\left[ \frac{\sqrt{B_n}}{\sqrt{B}}\, \frac{|\sqrt{C}- \sqrt{C_n}|}{\sqrt{C}}  \right]
= E\left[ \frac{\sqrt{B_n}}{\sqrt{B}} \right] \, E\left[ \frac{|\sqrt{C}- \sqrt{C_n}|}{\sqrt{C}} \right] \notag \\
&\leq& \frac{1}{12} \left\{ \frac{1}{132} \, \left[ \frac{5}{2} (C_3 +C_4) \right]^{\frac{1}{2}} + 1 \right\} 
\left[ \frac{5}{2} (C_3 +C_4) \right]^{\frac{1}{2}} \, \frac{1}{n}. \label{eq:difference2}  
\end{eqnarray}
By H\"{o}lder's inequality and Proposition \ref{prop:secondmoment},
\begin{eqnarray}
&&E\left[\frac{|A_n -A|}{\sqrt{B C}}\right]
\leq \left\{ E\left[ (A_n -A)^{2} \right] \right\}^{\frac{1}{2}} \, \left\{E\left[ B^{-1}C^{-1} \right] \right\}^{\frac{1}{2}} \notag  \\
&=& \left\{E\left[ (A_n -A)^{2} \right] \right\}^{\frac{1}{2}} \,  \left\{E[B^{-1}] E[C^{-1}]\right\}^{\frac{1}{2}}
\leq \left(\frac{5}{72 n^2}\right)^{\frac{1}{2}} \cdot C_1=   \frac{1}{6} \left(\frac{5}{2}\right)^{\frac{1}{2}} C_1 \,\frac{1}{n} . \label{eq:difference3} 
\end{eqnarray}
Combining \eqref{eq:differenceoffraction}, \eqref{eq:difference1}, \eqref{eq:difference2} and \eqref{eq:difference3} yields
\begin{equation*}
E\left[\left| \frac{A_n}{\sqrt{B_n C_n}} -  \frac{A}{\sqrt{B C}} \right|\right] 
\leq \frac{C_5}{n}.
\end{equation*}
This completes the proof.
\end{proof}

\section{Future work.} \label{sec6}
In this section, we present conjectures which, while beyond the scope of the present paper, should constitute opportunities for future work which should be tractable given some known tools and techniques in the analysis on Wiener chaos, and could have interesting applications to statistical testing based on paths of time series. 

The reader can refer to Section \ref{sec:coupling} for conjectures on convergence rates, and their implications, regarding the distinction between random walks and other types of time series. Those conjectures would apply to statistics which can be related to the Wasserstein distance.

Going beyond them, we conjecture that, for practical purposes, the convergence of $\theta_n$ to $\theta$ also occurs in total variation at the same rate as in Wasserstein distance, in the sense that the probability law of $\theta_{n}$ converges at the rate $r(n):= c n^{-1}$ for some constant $c$ 
though this may be harder to establish except empirically or via simulations.  The practical conjecture, that extends from the Wasserstein to the total variation distance, would be significant for several reasons, including because the total variation distance is an upper bound on the Kolmogorov distance, but only the square root of the Wasserstein distance bounds the Kolmogorov distance. As the latter is the supremum norm for the distance between cumulative distribution functions (CDFs), an application of the practical conjecture, using specifically the implication for the Kolmogorov distance, would be as follows. An upper bound of order of magnitude $10^{-2}$, say, could legitimately imply that any estimate on the $\alpha$-th percentile of $\theta$ could result in the same estimate on the $(\alpha -10^{-2})$-th percentile of $\theta_{n}$. One could thus build a test of independence of two (Gaussian) random walks of length $n$ where the rejection region at the confidence level $\alpha$ could be equated to the rejection region using the CDF of $\theta$ at the confidence level $\alpha + 10^{-2}$ as long as $r(n)<10^{-2}$. This argument could take into account the multiplicative constant $c$ in the speed of convergence $r(n)$, which could also be determined from simulations. Without our conjecture on total variation rate of convergence, using instead our Theorem \ref{thm:boundwassertain}, this strategy for rejection regions at level $\alpha + 10^{-2}$ would follow from $r(n)^{1/2}<10^{-2}$, because, as we mentioned, the Wasserstein distance only bounds the square root of the Kolmogorov distance.  

Other options for conjectures for statistical testing could include studying the speed of convergence of moment ratios of paths, such as a kurtosis-type statistic, and their fluctuations. Though this is also beyond the scope of this paper, we conjecture that, unlike the limit of the law of $\theta_{n}$ itself, whose numerator and denominator converge in the second chaos, the polarization of an empirical Kurtosis for two Gaussian random walks has normal fluctuations. Such a study could use either the technique presented in Section \ref{sec:convergenceinWasserstein} via bounding the negative moments of the denominator from its moment-generating function, or the so-called optimal fourth moment theorem (\cite{Nourdin}), where the speed of convergence of normal fluctuation for chaos sequences is known sharply in total variation.

We also suspect that the convergence phenomena we uncover here in the previous section are not restricted to Gaussian random walks, but hold for a wide range of random walks and other processes, including walks with other step distributions. Because of the heavy reliance on the Gaussian property in our work, particularly to be able to work in the second Wiener chaos, using non-Gaussian step distributions would require different tools. However, going from Gaussian random walks and Wiener processes to other Gaussian time series and their continuous limits could preserve a number of the tools we present here. For instance, we rely on the extraordinary convenience of Lemma \ref{lm:momentinverse} and the explicit nature of the corresponding moment-generating function, to estimate negative moments, but this can be done by other means for other Gaussian processes and their discrete-time observations. Similarly, as mentioned in Section \ref{sec:coupling}, we use the convenience of being able to calculate the exact value of the $L^2(\Omega)$ distance between the constituent elements of $\theta$ and $\theta_n$ (e.g. by employing Faulhaber's formula for the partial sum of the powers of integers). But these expressions can be estimated nearly as precisely, using the kernel representations, by invoking comparisons between series and Riemann integrals, with error terms of the same order as the second-order terms in Propositions \ref{prop:secondmoment} and \ref{prop:upperboundforBn-B}.

\newpage

\bibliographystyle{plain}
\bibliography{reference_Yule}

\newpage

\section{Appendix}

This appendix proves Lemma \ref{lem1} and Proposition \ref{prop:identicaldistribution}.

\subsection{Proof of Lemma \ref{lem1}}
For simplicity, we will start with $d_{n}(n^2 \lambda)$, for $n \geq 5$. From the definition of $d_{n}(\lambda)$, we have
\begin{eqnarray}
&& d_{n}(n^2 \lambda) = \det\left( I_{n-1} - n^2 \lambda K_{n}  \right)
= \det\left( I_{n-1} - \lambda\{ n \min(j,k) - jk \}_{j,k=1}^{n-1} \right) \notag  \\
&=& \left|
\begin{matrix}
 1- (n-1)\lambda,  & -(n-2) \lambda,  & -(n-3) \lambda, & -(n-4) \lambda, & \cdots,  & -\lambda  \\
 -(n-2)\lambda,  & 1-2(n-2) \lambda,  & -2(n-3) \lambda, & -2(n-4)\lambda, & \cdots,  & -2\lambda \\
 -(n-3)\lambda,  & -2(n-3) \lambda,  & 1-3(n-3) \lambda, & -3(n-4)\lambda, & \cdots, & -3\lambda \\
 -(n-4)\lambda,  & -2(n-4) \lambda,  & -3(n-4) \lambda, & 1-4(n-4)\lambda, & \cdots,  & -4\lambda \\
 \vdots &  \vdots & \vdots & \vdots & \ddots  &\vdots \\
 -\lambda,& -2\lambda, & -3\lambda, & -4\lambda, &\cdots, & 1-(n-1)\lambda 
\end{matrix}
\right| \label{eq:originalmatrix} \\  [1mm]
&=& \left|
\begin{matrix}
1- (n-1)\lambda,  & -(n-2) \lambda,  & -(n-3) \lambda, & -(n-4) \lambda, & \cdots,  & -\lambda  \\
n\lambda -2, & 1, & 0, & 0, &\cdots, & 0  \\
2n\lambda -3, & n\lambda , & 1, & 0, &\cdots, & 0  \\
3n\lambda -4, & 2 n \lambda, & n\lambda, & 1, &\cdots, & 0  \\
\vdots &  \vdots & \vdots & \vdots & \ddots  &\vdots  \\
(n-2)n\lambda -(n-1), & (n-3)n\lambda, & (n-4)n\lambda, & (n-5)n\lambda, &\cdots, & 1
\end{matrix}
\right| \label{eq:matrixtransform1}  \\ [1mm]
&=& \left|
\begin{matrix}
1- (n-1)\lambda,  & 0,  & 0,  & \cdots,  & -\lambda  \\
n\lambda -2, & 1, & 0,  &\cdots, & 0  \\
2n\lambda -3, & n\lambda , & 1,  &\cdots, & 0  \\
3n\lambda -4, & 2 n \lambda, & n\lambda,  &\cdots, & 0  \\
\vdots &  \vdots & \vdots &  \ddots  &\vdots  \\
(n-2)n\lambda -(n-1), & (n-3)n\lambda-(n-2), & (n-4)n\lambda-(n-3),  &\cdots, & 1
\end{matrix}
\right|, \label{eq:matrixtransform2}
\end{eqnarray}
where from \eqref{eq:originalmatrix} to \eqref{eq:matrixtransform1}, we add $(-j)\times$ the first row to the $j$-th row, for $j= 2,3,\cdots, n-1$. From \eqref{eq:matrixtransform1} to \eqref{eq:matrixtransform2}, we add $-(n-j)\times$ the last column to the $j$-th column, for $j=2,3\cdots,n-2$. We expand the determinant in \eqref{eq:matrixtransform2} by its first row and obtain
\begin{eqnarray}
&& d_{n}(n^2 \lambda) = 1-(n-1)\lambda + (-1)^{n+1} \lambda \notag \\
&& \cdot
\left|
\begin{matrix}
n\lambda -2, & 1, & 0,  &\cdots, & 0  \\
2n\lambda -3, & n\lambda , & 1,  &\cdots, & 0  \\
3n\lambda -4, & 2 n \lambda, & n\lambda,  &\cdots, & 0  \\
\vdots &  \vdots & \vdots &  \ddots  &\vdots  \\
(n-3)n\lambda -(n-2), & (n-4)n\lambda, & (n-5)n\lambda,  &\cdots, & 1  \\
(n-2)n\lambda -(n-1), & (n-3)n\lambda-(n-2), & (n-4)n\lambda-(n-3),  &\cdots, & n\lambda -2
\end{matrix}
\right|.  \label{eq:matrixtransform3}
\end{eqnarray}
Further, the determinant in \eqref{eq:matrixtransform3} is equal to
\begin{eqnarray}
  &=& \left|
  \begin{matrix}
   n\lambda -2, & 1, & 0, & \cdots, & 0, & 0, & 0   \\
   n\lambda -1, & n\lambda -1, & 1,  & \cdots, & 0, & 0, & 0   \\
   n\lambda -1, & n\lambda , & n\lambda -1,  & \cdots, & 0, & 0, & 0  \\
   \vdots & \vdots & \vdots & \ddots & \vdots & \vdots & \vdots   \\
   n\lambda-1, & n\lambda,  & n\lambda, & \cdots, & n\lambda-1, & 1, & 0  \\
   n\lambda-1, & n\lambda, & n\lambda, & \cdots, & n\lambda, & n\lambda-1, & 1  \\
   n\lambda-1, & n\lambda-(n-2), & n\lambda-(n-3), & \cdots, & n\lambda-4, & n\lambda-3, & n\lambda-3
  \end{matrix}
 \right|  \label{eq:matrixtransform4} \\  [1mm]
 &=& \left|
  \begin{matrix}
   n\lambda -2, & 1, & 0, & \cdots, & 0, & 0, & 0   \\
   1, & n\lambda -2, & 1,  & \cdots, & 0, & 0, & 0   \\
   0, & 1 , & n\lambda -2,  & \cdots, & 0, & 0, & 0   \\
   \vdots & \vdots & \vdots & \ddots & \vdots & \vdots & \vdots   \\
   0, & 0,  & 0, & \cdots, & n\lambda-2, & 1, & 0  \\
   0, & 0, & 0, & \cdots, & 1, & n\lambda-2, & 1  \\
   0, & -(n-2), & -(n-3), & \cdots, & -4, & -2, & n\lambda-4 
  \end{matrix}
 \right|.  \label{eq:matrixtransform5}
\end{eqnarray}
From \eqref{eq:matrixtransform3} to \eqref{eq:matrixtransform4}, we add $(-1)\times$ $(j-1)$-th row to $j$-th row for $j=n-2, n-3, \cdots, 2$. From \eqref{eq:matrixtransform4} to \eqref{eq:matrixtransform5}, similarly, we add $(-1)\times$ $(j-1)$-th row to $j$-th row for $j=n-2, n-3, \cdots, 2$.

Before proceeding to calculate $d_{n}(n^2 \lambda)$, we pause here to introduce a new determinant, closely related to $d_{n}(n^2 \lambda)$. Let us denote by $p_{n}(\lambda)$ the following $(n-2)\times (n-2)$ determinant:
\begin{equation}
\left|
  \begin{matrix}
   n\lambda -2, & 1, & 0, & \cdots, & 0, & 0, & 0   \\
   1, & n\lambda -2, & 1,  & \cdots, & 0, & 0, & 0   \\
   0, & 1 , & n\lambda -2,  & \cdots, & 0, & 0, & 0   \\
   \vdots & \vdots & \vdots & \ddots & \vdots & \vdots & \vdots  \\
   0, & 0,  & 0, & \cdots, & n\lambda-2, & 1, & 0  \\
   0, & 0, & 0, & \cdots, & 1, & n\lambda-2, & 1  \\
   -(n-1), & -(n-2), & -(n-3), & \cdots, & -4, & -2, & n\lambda-4 
  \end{matrix}
 \right|.  \label{eq:matrixdefforp}
\end{equation}
As mentioned, we introduce this determinant to compensate for the break in symmetry in \eqref{eq:matrixtransform5} because of the zero in the lower-left-hand corner there. We may easily verify that $p_{n}(\lambda) = \eqref{eq:matrixtransform5} + (-1)^{n}(n-1)$. Note that, from the expression of $d_n(n^2 \lambda)$ in \eqref{eq:matrixtransform3} and \eqref{eq:matrixtransform5}, we have
\begin{eqnarray}
d_{n}(n^2 \lambda) &=& 1-(n-1)\lambda + (-1)^{n+1} \lambda \times \eqref{eq:matrixtransform5}  \notag \\
&=& 1-(n-1)\lambda + (-1)^{n+1} \lambda \left( p_{n}(\lambda) - (-1)^{n}(n-1) \right) \notag \\
&=& 1 + (-1)^{n+1} \lambda \, p_{n}(\lambda). \label{eq:relationdandp}
\end{eqnarray}

By \eqref{eq:relationdandp}, the problem of calculating $d_{n}(n^2 \lambda)$ is converted to the problem of calculating $p_{n}(\lambda)$. To calculate $p_{n} (\lambda)$, our strategy is to derive a second-order recursion formula for $p_{n}(\lambda/n)$, see \eqref{eq:iteratedfunctionforp} below.
In what follows, we derive an explicit expression for $p_{n}(\lambda)$. 

For $n \geq 7$, we expand the determinant \eqref{eq:matrixdefforp} by its first column and obtain
\begin{eqnarray}
p_{n}(\lambda) &=& (n\lambda -2)\,p_{n-1}\left( \frac{n}{n-1} \lambda \right) + (-1)^{n}(n-1) \notag  \\
&& -1 \cdot
\left|
  \begin{matrix}
   1, & 0, & 0, & \cdots, & 0, & 0, & 0   \\
   1, & n\lambda -2, & 1,  & \cdots, & 0, & 0, & 0   \\
   0, & 1 , & n\lambda -2,  & \cdots, & 0, & 0, & 0   \\
   \vdots & \vdots & \vdots & \ddots & \vdots & \vdots & \vdots   \\
   0, & 0,  & 0, & \cdots, & n\lambda-2, & 1, & 0  \\
   0, & 0, & 0, & \cdots, & 1, & n\lambda-2, & 1  \\
   -(n-2), & -(n-3), & -(n-4), & \cdots, & -4, & -2, & n\lambda-4 
  \end{matrix}
 \right|.  \label{eq:matrixtransform6}
\end{eqnarray}
If we expand the determinant in \eqref{eq:matrixtransform6} by its first row, then it is exactly $p_{n-2}\left(\frac{n}{n-2} \lambda\right)$. Hence,
\begin{equation*}
p_{n}(\lambda) = (n\lambda -2)\, p_{n-1}\left( \frac{n}{n-1} \lambda \right) - p_{n-2}\left(  \frac{n}{n-2} \lambda \right)+(-1)^{n}(n-1).
\end{equation*}
In the above equation, we make a change of variables from $\lambda$ to $\lambda/n$ and obtain
\begin{equation*}
p_{n}\left(\frac{\lambda}{n}\right) = (\lambda -2)p_{n-1}\left(\frac{\lambda}{n-1}\right)
- p_{n-2}\left(\frac{\lambda}{n-2}\right) + (-1)^{n}(n-1).
\end{equation*}
For $\lambda \neq 0$, rearranging the above equation yields
\begin{eqnarray}
(-1)^{n}\, p_{n}\left(\frac{\lambda}{n}\right) - \frac{1}{\lambda}\,n
&=& -(\lambda -2) \left[ (-1)^{n-1} p_{n-1}\left(\frac{\lambda}{n-1}\right) - \frac{1}{\lambda}\, (n-1) \right]  \notag \\
&& - \left[ (-1)^{n-2} p_{n-2}\left(\frac{\lambda}{n-2}\right) - \frac{1}{\lambda}\,(n-2) \right]. \label{eq:iteratedfunctionforp}
\end{eqnarray}
The above iterative formula of $(-1)^{n} p_{n} (\lambda/n) - n/\lambda$ tells us that for $\lambda \neq 0$ or $4$, it must have the following form:
\begin{equation*}
C_1 \cdot \left( - \frac{(\lambda -2) + \sqrt{(\lambda-2)^2 - 4}}{2} \right)^{n}
+ C_2 \cdot \left( - \frac{(\lambda -2) - \sqrt{(\lambda-2)^2 - 4}}{2} \right)^{n},
\end{equation*}
where $C_1$ and $C_2$ are two constants.
Direct calculation gives
\begin{eqnarray*}
&& (-1)^{5}\, p_{5}\left(\frac{\lambda}{5}\right) - \frac{5}{\lambda} = -(\lambda -4)\left( (\lambda-2)^{2} +1 \right) - \frac{5}{\lambda},  \\
&& (-1)^{6}\, p_{6}\left(\frac{\lambda}{6}\right) - \frac{6}{\lambda} = (\lambda -4) (\lambda-2)^{3} +3  - \frac{6}{\lambda}.
\end{eqnarray*}
Then the constants $C_1$ and $C_2$ can be determined as
\begin{equation*}
C_1 = - C_2 = \frac{1}{\lambda \sqrt{(\lambda-2)^2 -4}}.
\end{equation*}
Hence, for $n \geq 5$ and $\lambda \neq 0$ or $4$,
\begin{eqnarray}
&& (-1)^{n}\, p_{n}\left(\frac{\lambda}{n}\right) - \frac{1}{\lambda}\,n  \notag \\
&=& \frac{1}{\lambda \sqrt{(\lambda-2)^2 -4}}
 \left[ \left( - \frac{(\lambda -2) + \sqrt{(\lambda-2)^2 - 4}}{2} \right)^{n} - \left( - \frac{(\lambda -2) - \sqrt{(\lambda-2)^2 - 4}}{2} \right)^{n} \right] \notag \\ [1mm]
&=& \frac{(-1)^n}{2^{n-1} \, \lambda} \, \sum_{k=1}^{\lceil n/2 \rceil} {n \choose 2k-1} \, (\lambda-2)^{n-(2k-1)} \, \left( (\lambda-2)^2 -4 \right)^{k-1}.  \label{eq:representationforscallingp}{}
\end{eqnarray}
Combining \eqref{eq:relationdandp} and \eqref{eq:representationforscallingp} yields that, for $n \geq 5$ and $\lambda \neq 0$ or $4n$, \eqref{eq:representationfordinmaindocument} holds.
Since both sides of \eqref{eq:representationfordinmaindocument} are continuous function of $\lambda$, it also holds for every $\lambda \in \mathds{R}$ and $n \geq 5$. It is straightforward to verify that \eqref{eq:representationfordinmaindocument} also holds for $n = 2,3$ and $4$, hence, \eqref{eq:representationfordinmaindocument} holds for all $n \geq 2$.

\subsection{Proof of Proposition \ref{prop:identicaldistribution}}
We first note that statement (c) follows immediately from statements (a) and (b). We first prove statement (b). By integration by parts for the Wiener integral, we have
\begin{eqnarray*}
\int_{0}^{t} M(s,t) \, dW_{1}(s) 
&=& \int_{0}^{t} (s-st) \, dW_{1}(s)  
= (t - t^{2})W_1(t) - (1-t) \int_{0}^{t} W_{1}(s) \, ds.
\end{eqnarray*}
Then
\begin{eqnarray}
&& \int_{0}^{1} \int_{0}^{t} M(s,t) \, dW_{1}(s) \, dW_{2}(t) \notag \\
&=& \int_{0}^{1} (t-t^2) W_1(t) \,dW_{2}(t) - \int_{0}^{1} (1-t) \left(\int_{0}^{t} W_{1}(s)\, ds \right) \,  dW_{2}(t). \label{eq:1}
\end{eqnarray}
Applying It\^{o}'s lemma to $(1-t)W_{2}(t) \int_{0}^{t} W_{1}(s)\, ds$ yields
\begin{eqnarray*}
0 &=& - \int_{0}^{1} \left(\int_{0}^{t} W_{1}(s) \, ds \right) W_{2}(t) \, dt 
+ \int_{0}^{1} (1-t) \left(\int_{0}^{t} W_{1}(s)\, ds \right) \,  dW_{2}(t)  \\
 && + \int_{0}^{1} (1-t) W_1(t) W_{2}(t) \, dt.
\end{eqnarray*}
Together with \eqref{eq:1}, we have
\begin{eqnarray*}
&& \int_{0}^{1} \int_{0}^{t} M(s,t) \, dW_{1}(s) \, dW_{2}(t)  \\
&=& \int_{0}^{1} (t-t^2) W_1(t) \,dW_{2}(t) + \int_{0}^{1} (1-t) W_1(t) W_{2}(t) \, dt
- \int_{0}^{1} \left(\int_{0}^{t} W_{1}(s) \, ds \right) W_{2}(t) \, dt .
\end{eqnarray*}
Similarly,
\begin{eqnarray*}
&& \int_{0}^{1} \int_{0}^{s} M(s,t) \, dW_{2}(t) \, dW_{1}(s)  \\
&=& \int_{0}^{1} (t-t^2) W_{2}(t) \,dW_{1}(t) + \int_{0}^{1} (1-t) W_{1}(t) W_{2}(t) \, dt
- \int_{0}^{1} \left(\int_{0}^{t} W_{2}(s) \, ds \right) W_{1}(t) \, dt .
\end{eqnarray*}
Then
\begin{eqnarray}
A &=& \int_{0}^{1} \int_{0}^{t} M(s,t) \, dW_{1}(s) \, dW_{2}(t)
+ \int_{0}^{1} \int_{0}^{s} M(s,t) \, dW_{2}(t) \, dW_{1}(s)  \notag \\
&=& \int_{0}^{1} (t-t^2) W_1(t) \,dW_{2}(t) + \int_{0}^{1} (t-t^2) W_{2}(t) \,dW_{1}(t) 
 +  \int_{0}^{1} (2-2t) W_{1}(t) W_{2}(t) \, dt  \notag \\
&& - \int_{0}^{1} \left(\int_{0}^{t} W_{1}(s) \, ds \right) W_{2}(t) \, dt 
- \int_{0}^{1} \left(\int_{0}^{t} W_{2}(s) \, ds \right) W_{1}(t) \, dt.  \label{eq:Aform}
\end{eqnarray}
Applying It\^{o}'s lemma to $(t -t^2)\, W_{1}(t)\, W_{2}(t)$ yields
\begin{equation}
0 = \int_{0}^{1} (1-2t) W_{1}(t) W_{2}(t) \, dt + \int_{0}^{1} (t-t^2) W_{2}(t) \,dW_{1}(t)  
+ \int_{0}^{1} (t-t^2) W_1(t) \,dW_{2}(t).  \label{eq:byproduct1}
\end{equation}
Note that
\begin{equation*}
\int_{0}^{1} \left(\int_{0}^{t} W_{1}(s) \, ds \right) W_{2}(t) \, dt 
= \int_{0}^{1} \left(\int_{s}^{1} W_{2}(t) \, dt \right) W_{1}(s) \, ds
= \int_{0}^{1} \left(\int_{t}^{1} W_{2}(s) \, ds \right) W_{1}(t) \, dt,
\end{equation*}
where the first equality follows by interchanging the order of integrals and the second equality follows by substituting $(s,t)$ for $(t,s)$. We then calculate
\begin{eqnarray}
&& \int_{0}^{1} \left(\int_{0}^{t} W_{1}(s) \, ds \right) W_{2}(t) \, dt 
+ \int_{0}^{1} \left(\int_{0}^{t} W_{2}(s) \, ds \right) W_{1}(t) \, dt    \notag \\
&=& \int_{0}^{1} \left(\int_{t}^{1} W_{2}(s) \, ds \right) W_{1}(t) \, dt
  + \int_{0}^{1} \left(\int_{0}^{t} W_{2}(s) \, ds \right) W_{1}(t) \, dt  \notag \\
&=& \int_{0}^{1} \left(\int_{0}^{1} W_{2}(s) \, ds \right) W_{1}(t) \, dt 
= \int_{0}^{1} W_{1}(t) \, dt   \int_{0}^{1} W_{2}(s) \, ds.  \label{eq:byproduct2}
\end{eqnarray}
Combining \eqref{eq:Aform}, \eqref{eq:byproduct1} and \eqref{eq:byproduct2}, \eqref{eq:A} follows. Since \eqref{eq:B} and \eqref{eq:C} are symmetric, we need only prove \eqref{eq:B}, and then \eqref{eq:C} will follow similarly. By a similar argument to that in the derivation of \eqref{eq:Aform}, we have
\begin{eqnarray}
B &=& 2 \int_{0}^{1} (t-t^2) W_1(t) \,dW_{1}(t)
 + \int_{0}^{1} (2-2t) \, W_{1}^{2}(t) \, dt   \notag \\
 && - 2 \int_{0}^{1} \left(\int_{0}^{t}  W_{1}(s) \, ds \right) W_{1}(t) \, dt
  + \int_{0}^{1} M(t,t)\, dt. \label{eq:Bform}
\end{eqnarray}
Applying It\^{o}'s lemma to $(t -t^2)\, W_{1}^{2}(t)$ yields
\begin{equation}
0 = \int_{0}^{1} (1-2t) \, W_{1}^{2}(t) \, dt
 +  2 \int_{0}^{1} (t-t^2) W_1(t) \,dW_{1}(t) + \int_{0}^{1} (t-t^2)\,dt.  \label{eq:byproduct3}
\end{equation}
Further,
\begin{eqnarray}
\int_{0}^{1} \left(\int_{0}^{t}  W_{1}(s) \, ds \right) W_{1}(t) \, dt
&=& \int_{0}^{1} \left(\int_{0}^{t}  W_{1}(s) \, ds \right) \,d\left(\int_{0}^{t}  W_{1}(s) \, ds \right)  \notag \\
&=& \frac{1}{2} \left( \int_{0}^{1}  W_{1}(s) \, ds \right)^2. \label{eq:byproduct4}
\end{eqnarray}
Combining \eqref{eq:Bform}, \eqref{eq:byproduct3} and \eqref{eq:byproduct4}, \eqref{eq:B} follows.

In the remainder of the proof, we will show that
\begin{eqnarray}
A_n  &=&  \sum_{j,k=1}^{n} M\left( \frac{j-1}{n}, \frac{k-1}{n} \right) \left( W_1\left( \frac{j}{n}\right) -  W_1\left( \frac{j-1}{n}\right) \right) \left( W_2\left( \frac{k}{n}\right) -  W_2\left( \frac{k-1}{n}\right) \right),  \label{eq:An}  \\
B_n  &=&  \sum_{j,k=1}^{n} M\left( \frac{j-1}{n}, \frac{k-1}{n} \right) \left( W_1\left( \frac{j}{n}\right) -  W_1\left( \frac{j-1}{n}\right) \right) \left( W_1\left( \frac{k}{n}\right) -  W_1\left( \frac{k-1}{n}\right) \right),  \label{eq:Bn}  \\
C_n  &=&  \sum_{j,k=1}^{n} M\left( \frac{j-1}{n}, \frac{k-1}{n} \right) \left( W_2\left( \frac{j}{n}\right) -  W_2\left( \frac{j-1}{n}\right) \right) \left( W_2\left( \frac{k}{n}\right) -  W_2\left( \frac{k-1}{n}\right) \right).  \label{eq:Cn}  
\end{eqnarray}
Note that $ W_1\left( \frac{j}{n}\right) -  W_1\left( \frac{j-1}{n}\right)$, $W_2\left( \frac{k}{n}\right) -  W_2\left( \frac{k-1}{n}\right)$, $j,k = 1,2 \dots, n$ are mutually independent Gaussian random variables with distribution $\mathcal{N}(0, 1/n)$. By Section \ref{sec:jointmgf} (see line \eqref{eq:quadraticformcoefficients}), it follows easily that $Z_{12}^{n}/n$, $Z_{11}^{n}/n$ and $Z_{22}^{n}/n$ are quadratic forms of the random variables $(X_{1}/\sqrt{n}, X_{2}/\sqrt{n}, \dots, X_{n}/\sqrt{n})$ and $(Y_{1}/\sqrt{n}, Y_{2}/\sqrt{n}, \dots, Y_{n}/\sqrt{n})$ with same coefficients as $A_n$, $B_n$ and $C_n$ respectively. Thus, statement (a) of Proposition \ref{prop:identicaldistribution} follows immediately.

For simplicity, let $u_j$ denote $j/n$ for $j=0,1,2, \dots, n$. We proceed to calculate
\begin{eqnarray}
&&\int_{0}^{1} \int_{0}^{t} M_{n}(s,t) \, dW_{1}(s) \, dW_{2}(t) \notag \\ 
&=& \int_{0}^{1} \int_{0}^{t} \sum_{ j,k =1}^{n} M\left(u_{j-1}, u_{k-1}\right) \mathds{1}_{\{ u_{j-1} < s \leq u_{j}\}} 
\mathds{1}_{\{ u_{k-1} < t \leq u_{k}\}} \, dW_{1}(s) \, dW_{2}(t) \notag  \\
&=& \sum_{ j<k } \int_{0}^{1} \int_{0}^{t} M\left(u_{j-1}, u_{k-1}\right) \mathds{1}_{\{ u_{j-1} < s \leq u_{j}\}} 
\mathds{1}_{\{ u_{k-1} < t \leq u_{k}\}} \, dW_{1}(s) \, dW_{2}(t)  \notag \\
&& + \sum_{j=1}^{n} \int_{0}^{1} \int_{0}^{t} M\left(u_{j-1}, u_{j-1}\right) \mathds{1}_{\{ u_{j-1} < s \leq u_{j}\}} 
\mathds{1}_{\{ u_{j-1} < t \leq u_{j}\}} \, dW_{1}(s) \, dW_{2}(t),  \label{eq:sumation}
\end{eqnarray}
where the equality holds because the term with indices satisfying $j>k$ is $0$. The first term on the right-hand side of \eqref{eq:sumation} is
\begin{eqnarray}
&& \sum_{j<k} \int_{0}^{1} M\left(u_{j-1}, u_{k-1}\right) \left( W_1\left( u_{j}\right) -  W_1\left( u_{j-1}\right) \right) \mathds{1}_{\{ u_{k-1} < t \leq u_{k} \}}  \, dW_{2}(t) \notag \\
&=& \sum_{j<k} M\left(u_{j-1}, u_{k-1}\right) \left( W_1\left( u_{j}\right) -  W_1\left( u_{j-1}\right) \right) \left( W_2\left( u_{k} \right) -  W_2\left( u_{k-1} \right) \right). \label{eq:firstterm}
\end{eqnarray}
The second term on the right-hand side of \eqref{eq:sumation} is
\begin{eqnarray}
&& \sum_{j=1}^{n} \int_{0}^{1} M\left( u_{j-1}, u_{j-1} \right) \left(W_{1}(t) - W_{1}\left(u_{j-1}\right)\right) \mathds{1}_{\{ u_{j-1} < t \leq u_{j}\}} \, dW_{2}(t) \notag \\
&=& \sum_{j=1}^{n} M\left( u_{j-1}, u_{j-1} \right) \left( \int_{u_{j-1}}^{u_{j}} W_{1}(t) \,dW_2(t) - W_{1}\left(u_{j-1}\right) \left( W_{2}\left(u_{j}\right) -W_{2}\left(u_{j-1}\right) \right) \right).
 \label{eq:secondterm}
\end{eqnarray}
Combining \eqref{eq:sumation}, \eqref{eq:firstterm} and \eqref{eq:secondterm} yields
\begin{eqnarray}
&&\int_{0}^{1} \int_{0}^{t} M_{n}(s,t) \, dW_{1}(s) \, dW_{2}(t) \notag \\
&=& \sum_{j<k} M\left(u_{j-1}, u_{k-1}\right) \left( W_1\left( u_{j}\right) -  W_1\left( u_{j-1}\right) \right) \left( W_2\left( u_{k} \right) -  W_2\left( u_{k-1} \right) \right) \notag \\
&& + \sum_{j=1}^{n} M\left( u_{j-1}, u_{j-1} \right) \left( \int_{u_{j-1}}^{u_{j}} W_{1}(t) \,dW_2(t) - W_{1}\left(u_{j-1}\right) \left( W_{2}\left(u_{j}\right) -W_{2}\left(u_{j-1}\right) \right) \right).
\label{eq:Anpart1}
\end{eqnarray}
Similarly,
\begin{eqnarray}
&&\int_{0}^{1} \int_{0}^{s} M_{n}(s,t) \, dW_{2}(t) \, dW_{1}(s) \notag \\
&=& \sum_{j<k} M\left(u_{j-1}, u_{k-1}\right) \left( W_2\left( u_{j}\right) -  W_2\left( u_{j-1}\right) \right) \left( W_1\left( u_{k} \right) -  W_1\left( u_{k-1} \right) \right) \notag \\
&& + \sum_{j=1}^{n} M\left( u_{j-1}, u_{j-1} \right) \left( \int_{u_{j-1}}^{u_{j}} W_{2}(t) \,dW_1(t) - W_{2}\left(u_{j-1}\right) \left( W_{1}\left(u_{j}\right) -W_{1}\left(u_{j-1}\right) \right) \right).
\label{eq:Anpart2}
\end{eqnarray}
Combining \eqref{eq:Anpart1} and \eqref{eq:Anpart2} and rearranging terms gives
\begin{eqnarray}
 A_n &=&  
 \sum_{j,k=1}^{n} M\left(u_{j-1}, u_{k-1}\right) \left( W_2\left( u_{j}\right) -  W_2\left( u_{j-1}\right) \right) \left( W_1\left( u_{k} \right) -  W_1\left( u_{k-1} \right) \right)  \notag \\
&& +  \sum_{j=1}^{n} M\left( u_{j-1} , u_{j-1} \right) \Bigg[ \int_{u_{j-1}}^{u_{j}} W_{1}(t) \,dW_2(t) + \int_{u_{j-1}}^{u_{j}} W_{2}(t) \,dW_1(t)\notag \\
&& \quad \quad -  \left( W_1\left( u_{j}\right)  W_2\left( u_{j}\right) - W_1\left( u_{j-1}\right)  W_2\left( u_{j-1}\right)\right)  \Bigg].  \label{eq:Anform}
\end{eqnarray}
Applying It\^{o}'s lemma to $W_{1}(t) W_{2}(t)$ yields
\begin{eqnarray*}
&&  W_1\left( u_{j}\right)  W_2\left( u_{j}\right) - W_1\left( u_{j-1} \right)  W_2\left( u_{j-1} \right) 
=   \int_{u_{j-1}}^{u_{j}} W_{2}(t) \,dW_1(t) + \int_{u_{j-1}}^{u_{j}} W_{1}(t) \,dW_2(t).
\end{eqnarray*}
Together with \eqref{eq:Anform}, \eqref{eq:An} follows. Since \eqref{eq:Bn} and \eqref{eq:Cn} are symmetric, we need only prove \eqref{eq:Bn}, and then \eqref{eq:Cn} will follow similarly. By a similar argument to that of the derivation of \eqref{eq:Anform}, $B_n$ equals
\begin{eqnarray*}
&& \sum_{j,k=1}^{n} M\left( u_{j-1}, u_{k-1} \right) \left( W_1\left(u_{j}\right) -  W_1\left( u_{j-1}\right) \right) \left( W_1\left( u_{k}\right) -  W_1\left( u_{k-1}\right) \right)  \notag \\
&& +  \sum_{j=1}^{n} M\left( u_{j-1}, u_{j-1} \right) \Bigg[ 2\int_{u_{j-1}}^{u_{j}} W_{1}(t) \,dW_1(t) 
  - \left( W_1^2\left( u_{j}\right)   - W_1^2\left( u_{j-1}\right)  \right) + \frac{1}{n}  \Bigg]. 
\end{eqnarray*}
Note that by applying It\^{o}'s lemma to $W_{1}^{2}(t)$
\begin{eqnarray*}
 W_1^2\left( u_{j}\right)  - W_1^2\left( u_{j-1}\right) 
 = 2\int_{u_{j-1}}^{u_{j}} W_{1}(t) \,dW_1(t)  + \frac{1}{n},
\end{eqnarray*}
from which \eqref{eq:Bn} immediately follows.

\end{document}